\documentclass[preprint, 11pt, nonatbib]{elsarticle}

\makeatletter
\@ifundefined{c@author}{}{\let\c@author\relax}
\makeatother

\makeatletter
\def\ps@pprintTitle{%
  \let\@oddfoot\@empty
  \let\@evenfoot\@empty
}
\makeatother

\usepackage[
backend=biber,
style=apa,
hyperref = true
]{biblatex}
\usepackage{setspace}
\usepackage{geometry}
\usepackage[T1]{fontenc}
\usepackage[colorlinks=true,
            linkcolor=blue,
            citecolor=blue,
            urlcolor=blue]{hyperref}

\usepackage{comment}
\usepackage{amsmath, amsfonts, nicefrac, booktabs, microtype, graphicx}
\usepackage{url, rotating, lscape, caption, subcaption, lipsum, threeparttable}
\usepackage{tikz}
\usepackage{pgfplots}
\usepgfplotslibrary{statistics}
\pgfplotsset{compat=1.18}

\usepackage{csvsimple}
\usetikzlibrary{matrix}
\usepackage[table]{xcolor}
\usepackage{algorithm2e}
\usepackage{titlesec}

\definecolor{CARcolor}{HTML}{C2C9D8}
\definecolor{ENTcolor}{HTML}{E9E0C9}
\definecolor{ENTcolordark}{HTML}{CCB782}
\definecolor{whiteback}{HTML}{FFFFFF}
\definecolor{mygreen}{HTML}{229954}
\definecolor{CAR1color}{HTML}{32107a}
\definecolor{ENT1color}{HTML}{81b5aa}
\definecolor{PLFcolor}{HTML}{81b5aa}
\definecolor{FNNcolor}{HTML}{a8b2c8}
\definecolor{SBMcolor}{HTML}{d9c9a3}

\titleformat{\section}
  {\normalfont\large\bfseries}
  {\thesection}{1em}{}

\titleformat{\subsection}
  {\normalfont\normalsize\bfseries}
  {\thesubsection}{1em}{}

\titleformat{\subsubsection}
  {\normalfont\normalsize\itshape}
  {\thesubsubsection}{1em}{}

\titleformat{\paragraph}[runin]{\normalfont\normalsize\bfseries}{\theparagraph}{1em}{}

\begin{document}

\begin{frontmatter}

\title{Handling Overtime Constraints in Mixed Integer Linear Programming for Surgical Scheduling: A Comparison of Neural Network and Classical Linearization Techniques}

\author[inst1,inst2]{Cindy Pistorius\corref{cor1}}
\cortext[cor1]{Corresponding author}
\ead{c.pistorius@erasmusmc.nl}

\author[inst2]{J. Theresia van Essen}
\ead{j.t.vanessen@tudelft.nl}

\address[inst1]{Erasmus University Medical Center, PO Box 2040, 3000 CA, Rotterdam, the Netherlands}
\address[inst2]{Delft Institute of Applied Mathematics, Delft University of Technology, Mekelweg 4, 2628 CD, Delft, the Netherlands}


\begin{abstract}
Uncertainty in surgery durations continues to be difficult to account for in operating room scheduling.
In particular, it remains complex to accurately incorporate uncertainty in surgical overtime constraints within mixed-integer linear programming (MILP) models.
Therefore, we propose a method that integrates feedforward neural networks (FNNs) into MILP models to approximate the total surgery duration in these overtime constraints.
The proposed approach is evaluated using real-life hospital data and compared against two classical approaches: scenario-based modelling and piecewise linear function approximations.
We demonstrate that with a relatively small FNN, we achieve competitive operating room schedules in terms of both solution quality and computational performance.
The FNN-based approach is the most computationally efficient with an optimality gap lower than 2\% in all cases, achieves the highest operating room utilization in six out of eight considered cases, and on average produces simulated overtime probabilities closest to the predefined target.\end{abstract}


\end{frontmatter}

\renewcommand{\arraystretch}{1.2}

\section{Introduction} \label{Sec:intro}
\noindent
One of the most expensive resources of a hospital are the operating rooms (ORs).
Improperly coordinated ORs can lead to substantial financial losses, as well as a disruption in patient treatment and dissatisfied personnel \parencite{Riahi2023ImprovingDuration, Rothstein2018OperatingEfficiency}. 
One key element of coordinating the ORs is the scheduling of surgeries. To create an OR schedule, many factors need to be taken into account, such as different departments and uncertainties of surgical processes. 
These processes can lead to various challenges on an operational, tactical and strategic level related to both elective and non-elective surgical patients \parencite{Rahimi2021AScheduling}. 
In particular, the operational level \parencite{Rahimi2021AScheduling, Abdelrasol2013AProblems} is affected, since it involves day-to-day scheduling decisions regarding the OR schedule.
This level covers the planning of surgical cases within the allocated time slots, considering the specific requirements of individual procedures. It addresses real-time resource allocation, along with managing any unforeseen changes such as emergencies or cancellations. 
The operational level is divided into two sub-levels: offline and online. The offline level covers the decisions made before the execution of tasks, such as the planning of surgical cases one week in advance. The online level considers (near) real-time decisions, such as emergencies and cancellations, while taking the current dynamic conditions and uncertainties into account.

In this paper, we focus on the scheduling of elective surgeries on the offline operational level. 
Within research, different techniques have been analysed to create surgical schedules. 
Mathematical models are among the most commonly used methodologies for these types of scheduling problems \parencite{Wang2021OperatingResearch}. 
One of these methodologies is Mixed Integer Linear Programming (MILP), which is a mathematical modelling approach that optimizes a linear objective function subject to a set of linear constraints, where a set of decision variables are restricted to be integers. 
MILP approaches remain popular due to their capacity to represent complex, constraint-rich scheduling problems while delivering reliable and interpretable results. They are also able to balance different objectives, which makes it possible to consider different criteria of stakeholders within one framework.
MILP models for scheduling surgeries are applicable to a range of objectives, such as optimizing the resource utilization to minimize patients' waiting times for surgery \parencite{Hachicha2018Two-MILPFacility},  
focusing on specific OR utilization rates \parencite{Clavel2016OperationDepartment} or minimizing costs while maximizing patient preferences \parencite{Ahmed2020ModelingProblem}. 
One of the challenges in scheduling surgical cases using an MILP model is accounting for the uncertainty of surgery durations.
This element is critical to the scheduling process, as surgical duration significantly impacts planned surgeries and affects not only the OR complex but also downstream departments.
Within the MILP framework, surgery durations are typically incorporated through the inclusion of overtime constraints, aimed at mitigating the risks of both underutilization and overtime across one or more ORs.
However, due to the linear formulation requirements of MILP constraints, it is challenging to incorporate the stochastic nature of these durations.

In research, many different strategies have been applied to consider this uncertainty inside MILP models. One approach is to use a pre-defined fixed value for each surgery type, such as the mean \parencite{Clavel2016OperationDepartment, Hachicha2018Two-MILPFacility, Jittamai2011ANo-shows}. This leads to an easy-to-implement set of constraints, where the total surgery duration for a specific OR is incorporated in overtime constraints as the sum of the assigned fixed values of the separate surgeries. The downside, however, is that these fixed values fail to capture the real-world variability of surgery durations, a crucial aspect of designing effective schedules.

To better account for this stochastic nature, other studies have used scenario-based modelling. With scenario-based modelling, the uncertainty is taken into account with a finite set of possibilities, i.e. scenarios.  
\textcite{Jebali2017AConstraints} and \textcite{Kamran2019AdaptiveApproach} both simulate different scenarios for the surgery duration based on the corresponding lognormal distribution.
With these scenarios, they minimize a certain objective that combines the outcomes across all scenarios by summing the scenario-specific objectives, each belonging to the same decision, to minimize the total cost across the scenarios.
\textcite{Rachuba2014AObjectives} use scenarios with equal probabilities of occurrence for the surgery duration to incorporate the characteristics of the lognormal distribution. Furthermore, they introduce a set of variables to account for scenario-dependent deviations in their problem. 
With scenario-based modelling, the uncertainty of stochastic variables can, thus, be addressed.
However, capturing the variability of surgery durations requires numerous scenarios, which can result in a high computational burden when many different surgery types need to be taken into account. 
\textcite{Doneda2024AScenarios} use pre-defined categories to create a limited number of total scenarios for the arrival and duration of emergency cases. This could be one strategy to address the computational burden. However, the balance between simplification and computational burden is challenging and can affect both the feasibility and performance of the model.

Another approach is to linearize the overtime constraint by approximating the total surgery duration for each OR with a normal distribution and apply a piecewise linear function. \textcite{Schneider2020SchedulingResources} have applied this technique to linearize their overtime constraints, which state that the probability of the total scheduled surgery duration in a specific OR being less than or equal to its capacity must be at least a pre-specified threshold.
They assume that the surgery duration is normally distributed, which results in a normally distributed total surgery duration for each OR within the MILP. 
This makes it possible  to apply the piecewise linear function and linearize the overtime constraint. A strength of this approach is its ability to incorporate the joint distribution of surgeries for each OR.
However, the lognormal distribution, identified as the most suitable distribution for surgery durations \parencite{Strum2000ModelingTimes}, cannot be used within this framework, since certain components of the lognormal distribution cannot be linearized with a piecewise linear function.

An alternative method to integrate non-linear constraints is to apply machine learning. With this approach, machine learning is used to learn a constraint or approximate known constraints based on available data \parencite{Fajemisin2024OptimizationSurvey}, after which the learned or approximated constraint is integrated in an optimization problem. 
This technique has been successfully applied in different case studies outside healthcare schedule optimization such as scheduling for the World Food Programme where palatability constraints are incorporated \parencite{Maragno2025Mixed-IntegerLearning}, chemotherapy optimization with predicted survival and toxicity levels \parencite{Maragno2025Mixed-IntegerLearning}, or auction optimization for which regression models are applied to learn revenue constraints \parencite{Verwer2017AuctionPrograms}.
The framework of \textcite{Fajemisin2024OptimizationSurvey} illustrates the versatile range of utilizing machine learning techniques to solve optimization problems. 
These techniques can also be extended to MILP models.
Because the total surgery duration is a continuous variable, in our research, only regression models are suitable for accounting for uncertainty within overtime constraints.
An additional important component is the approximated constraint's compatibility with the requirements of MILP models. 
\textcite{Maragno2025Mixed-IntegerLearning} state that neural networks can be represented with binary variables and big-M constraints, as long as the rectified linear unit (ReLU) activation function is used (or a linear activation function, as this function already has the linear properties). 
Other activation functions can be applied provided that the used solver supports these functions. 
Furthermore, they show that decision trees can be captured by linear equations, where each path in the tree can be represented using a set of constraints and binary variables. 
However, the number of constraints grows with the size of the tree. Random forests and ensemble trees can be expressed in the same manner, where a set of constraints is drawn up for each forest or ensemble as well as a set of global constraints.
These regression models make it possible to linearize the overtime constraints and to assume that the surgery duration follows a lognormal duration.
Under this assumption, the total surgery duration for each OR can be approximated by another lognormal distribution using the Fenton-Wilkinson method \parencite{Fenton1960TheSystems}.
This approach offers interesting possibilities, because the lognormal distribution captures the right-skewed nature of surgery durations more accurately, reducing the likelihood of underestimating the duration of longer procedures and thereby improving the reliability of scheduling. However, one drawback is that inaccurate model predictions may result in inefficient scheduling or even fail to satisfy the overtime constraints. Therefore, high performance of a predictive model in this context is crucial, and enough training data must be available to properly fit a regression model to approximate the lognormal distribution of the total surgery duration.

In this research, we propose a method that uses a neural network to formulate a linear approximation of the overtime constraints.
We compare this approach with two methods from literature, namely scenario-based modelling and using a piecewise linear approximation.
All three methods are compared in terms of solution quality and computational efficiency. 
The remainder of this paper is structured as follows: in Section \ref{Sec:problemdescr}, we describe the chosen elective scheduling problem at hand, which is derived from the elective surgery scheduling procedure at the Erasmus MC, an academic medical centre in the Netherlands. Section \ref{Sec:methods} covers the three different methods to take the overtime constraints into account within an MILP model. These methods are compared on their performance in Section \ref{Sec:results}.
This comparison is performed using surgical data from the Erasmus spanning 2.5 years.
The processing of the data is explained in Section \ref{sec:data_processing}.
We discuss our findings and reflect on the different approaches in Section \ref{Sec:discussion}.

\section{Problem description} \label{Sec:problemdescr}
\noindent
We begin by outlining the key assumptions for our model and presenting a generic MILP formulation for elective surgery scheduling.
We consider the situation where a hospital has dedicated ORs for each available specialty. Therefore, we formulate a generic MILP for a single specialty, which can be extended to multiple specialties. 
Furthermore, we assume that there are separate dedicated ORs for emergency surgeries, and therefore, do not consider non-elective surgeries.
Additionally, we assume that resources related to instruments and personnel are available and surgeries can be scheduled without any breaks in between.
The objective is to maximize OR utilization, while taking the capacity of each OR into account.

To start, we define $\mathcal{D} = \{0, 1,\dots,
D\}$ as the set of days in the planning horizon. 
For each specialty, a set of elective surgeries is given by set $\mathcal{S}$, where each surgery $s \in \mathcal{S}$ has a release date $r_s \in \mathbb{N}_0$, which is the earliest allowable surgery date. If surgery $s \in \mathcal{S}$ has a release date before the planning horizon, $r_s$ is set to 0.
Additionally, each surgery $s \in \mathcal{S}$ has a due date $q_s \in \mathbb{N}_0$, which is defined relative to the start of the planning horizon. 
As the planning horizon is indexed from day 0, a surgery $s \in \mathcal{S}$ due on the $n$-th calendar day has due date $q_s = n-1$. 
If a surgery does not have a due date, we set $q_s$ equal to the maximum due date among all surgeries plus 1.
For each surgery, we assign a binary parameter $p_s$, which is $1$ when surgery $s \in \mathcal{S}$ has a due date $q_s$ within the planning horizon $\mathcal{D}$ and $0$ otherwise. 
In addition, the surgery duration $\tilde{w}_s$ (in minutes) of surgery $s \in \mathcal{S}$ is a stochastic variable for which the mean surgery duration and accompanying variance are known, based on the corresponding surgery type.
Furthermore, the available set of ORs is denoted by $\mathcal{O}$. The maximum capacity in minutes of OR $o \in \mathcal{O}$ on day $d \in \mathcal{D}$ is then given by $C_{od}$.

Next, we introduce binary decision variables $X_{sod}$, which are $1$ when surgery $s \in \mathcal{S}$ is scheduled in operating room $o \in \mathcal{O}$ on day $d \in \mathcal{D}$ and $0$ otherwise. Each surgery can be scheduled at most once, which leads to the constraints:
\begin{equation}\label{eq:max_one}
    \sum_{o \in \mathcal{O}}\sum_{d \in \mathcal{D}} X_{sod} \leq 1, \quad \forall s \in \mathcal{S} .
\end{equation} 
A surgery cannot be scheduled before its release date. If a surgery has a due date within the planning horizon, the surgery has to be scheduled between the release and due date. This leads to the following constraints:
\begin{equation}
    \sum_{o \in \mathcal{O}}\sum_{d \,= \,0}^{r_s-1}X_{sod}=0, \quad \forall s \in \mathcal{S}, 
\end{equation}
\begin{equation}
    \sum_{o \in \mathcal{O}} \sum_{d\,=\,r_s}^{q_s} X_{sod} \geq p_s, \quad \forall s \in \mathcal{S}. 
\end{equation}


Furthermore, schedulers are generally not allowed to schedule a total mean surgery duration in a specific OR that exceeds its capacity. Therefore, we introduce the following capacity constraints, where $\bar{w}_s$ defines the sample mean of the surgery duration for each surgery $s \in \mathcal{S}$:
\begin{equation} \label{eq:overtime_samplemean}
    \sum_{s \in \mathcal{S}}X_{sod} \bar{w}_s \leq C_{od}, \quad \forall  o \in \mathcal{O}, \quad d \in \mathcal{D}.
\end{equation}
With these constraints, we do not consider surgeries that have a sample mean surgery duration $\bar{w}_s$ larger than the maximum available OR capacity: $\max_{o \in \mathcal{O}, d \in \mathcal{D}}\{C_{od}\}$. We assume that these rare cases are reviewed and scheduled manually by schedulers, and are therefore excluded.

Even though constraints \eqref{eq:overtime_samplemean} ensure that the total sample mean scheduled surgery duration may not exceed the capacity of an operating room, these capacity constraints do not guarantee that no overtime will occur, due to the stochastic nature of surgery durations.
Therefore, we introduce generic overtime constraints. Let $\tilde{k}_{od}$ be a stochastic variable that represents the total surgery duration scheduled at OR $o \in \mathcal{O}$ on day $d \in \mathcal{D}$. 
This leads to the following overtime constraints, which state that the probability on overtime should be less than or equal to a certain overtime probability $\alpha$ with $0 \leq \alpha \leq 1$, which can also be formulated as:
\begin{equation} \label{eq: generalovertime}
    P(\tilde{k}_{od} \leq C_{od}) \geq 1 - \alpha, \quad \forall  o \in \mathcal{O}, \quad d \in \mathcal{D}.
\end{equation}
At this stage, we have not explicitly specified $\tilde{k}_{od}$. Based on the solution method, the total surgery duration is further defined.

Lastly, we define the objective function. The goal of the model is to maximize OR utilization. We link this to the overall scheduled surgery time, since a total maximization of OR utilization leads to a maximization of total scheduled surgery duration. 
To ensure that we prioritize surgeries with an earlier due date when surgeries have equal surgery duration, we introduce a priority weight. Each surgery $s \in \mathcal{S}$ is assigned a priority $\frac{1}{q_s+1}$. This encourages the scheduling of surgeries with earlier due dates when surgery durations are otherwise equal. 
This results in the following objective function:
\begin{equation}\label{eq:generalobjfunc}
    \text{max}  \sum_{s \in \mathcal{S}}\sum_{o \in \mathcal{O}}\sum_{d \in \mathcal{D}}X_{sod} \left( \bar{w}_s+ \frac{1}{q_s+1} \right),
\end{equation}
where the sample mean $\bar{w}_s$ of the surgery duration for each scheduled surgery $s \in \mathcal{S}$ is included.
We do not take the stochastic behaviour of the surgery duration into account within the objective function to be able to compare the different methods based on their objective function value.

\section{Solution methods}\label{Sec:methods}
\noindent
In this section, we define three methods to incorporate the generic overtime constraints in the MILP model, formulated in Section \ref{Sec:problemdescr}. The first method we present is the linearization of the overtime constraints using a neural network, where we assume that the total surgery duration is lognormally distributed. 
Next, we discuss the linearization of the overtime constraints using a piecewise linear function, for which we assume that the total surgery duration is normally distributed. 
Lastly, we define the overtime constraints using a pre-defined set of scenarios, which we construct under the assumption that the surgeries are lognormally distributed. 

\subsection{Constraint linearization with neural network}
\noindent
The first method we elaborate on is the linearization of the overtime constraints by applying a machine learning algorithm.
This algorithm is used to predict the $(1-\alpha)$-percentile of the total surgery duration distribution.
For this method, we assume that the surgery durations are independent and identically distributed (i.i.d.), following a lognormal distribution, i.e. $\tilde{w}_{s} \sim \text{LogNormal}(\mu_{s}, \sigma^2_{s})$ for all $s \in \mathcal{S}$. 
To approximate the total surgery duration of surgeries scheduled in OR $o \in \mathcal{O}$ on day $d \in \mathcal{D}$, denoted as $\tilde{k}_{od}$, we need the approximation of the sum of the lognormal distributions of the individual surgeries. 
The Fenton-Wilkinson method \parencite{Fenton1960TheSystems} provides an effective way to approximate the sum of lognormal variables with another lognormal distribution.
The method approximates the distribution by matching the first two moments, i.e. the mean and the variance, with those of a lognormal distribution. 
This allows the estimation of the overall lognormal distribution of the total surgery duration for each OR $o \in \mathcal{O}$ and day $d \in \mathcal{D}$ based on the lognormal distributions of the individual surgeries.
Although this approach does not capture the exact distributional shape of the sum of the individual lognormal variables, it offers a computationally efficient approximation.

With the Fenton-Wilkinson method and the assumption that the surgery durations are i.i.d., we assume that $\tilde{k}_{od}$ is lognormally distributed in overtime constraints \eqref{eq: generalovertime}, i.e. $\tilde{k}_{od} \sim \text{LogNormal}(\mu_{od}, \sigma^2_{od})$ for all $o \in \mathcal{O}$ and $d \in \mathcal{D}$. We, then, introduce the following expressions:
\begin{equation} \label{eq:totmean}
    \mathop{\mathbb{E}}[\tilde{k}_{od}] =  \sum_{s \in \mathcal{S}} X_{sod}\mathop{\mathbb{E}}[\tilde{w}_s], \quad \forall  o \in \mathcal{O}, \quad d \in \mathcal{D},
\end{equation}
\begin{equation} \label{eq:totvar}
    \text{Var}(\tilde{k}_{od}) = \sum_{s \in \mathcal{S}} X_{sod} \text{Var}(\tilde{w}_s),  \quad \forall o \in \mathcal{O}, \quad d \in \mathcal{D}.
\end{equation}
The variables $\mathop{\mathbb{E}}[\tilde{w}_s]$ and $\text{Var}(\tilde{w}_s)$ are obtained as follows, where $\mu_s$ and $\sigma^2_s$ represent the parameters belonging to the lognormal distribution of the individual surgeries $\tilde{w}_s$:
\begin{equation} \label{eq:exp_var_lognormal}
\mathop{\mathbb{E}}[\tilde{w}_s] = e^{\mu_s + \frac{1}{2}\sigma^2_s} \quad \text{and} \quad \text{Var}(\tilde{w}_s)= \left( e^{\sigma^2_s} -1 \right) \left(\mathop{\mathbb{E}}[\tilde{w}_s]\right)^2.
\end{equation}
Following the relationships of the lognormal distribution, we find the lognormal parameters $\mu_{od}$ and $\sigma^2_{od}$ for OR $o \in \mathcal{O}$ and day $d \in \mathcal{D}$ using \eqref{eq:totmean} and \eqref{eq:totvar}:
\begin{equation}\label{eq:mu-sigma}
    \mu_{od} = \ln(\mathop{\mathbb{E}}[\tilde{k}_{od}]) - \frac{\sigma^2_{od}}{2}
\quad \text{and} \quad
    \sigma^2_{od} = \ln \left( \frac{\text{Var}(\tilde{k}_{od})}{\mathop{\mathbb{E}}[\tilde{k}_{od}]^2}+1 \right).
\end{equation}
With equations \eqref{eq:totmean}-\eqref{eq:mu-sigma}, we rewrite overtime constraints \eqref{eq: generalovertime} as follows:
\begin{equation} \label{eq:lognorm}
    P(\tilde{k}_{od} \leq C_{od}) = \Phi \left( \frac{\ln(C_{od})-\mu_{od}}{\sigma_{od}} \right) \geq 1 - \alpha, \quad \forall  o \in \mathcal{O}, \quad d \in \mathcal{D}.
\end{equation}
This can be expressed as
\begin{equation} \label{eq:overtimelog}
     e^{\mu_{od} + \sigma_{od}\Phi^{-1}(1-\alpha) } \leq C_{od},  \quad \forall  o \in \mathcal{O}, \quad d \in \mathcal{D},
\end{equation}
where $\Phi^{-1}(1-\alpha)$ indicates the $(1-\alpha)$-percentile of the normal distribution. 
This corresponds to the value that, when applied to the lognormal distribution with parameters $\mu_{od}$ and $\sigma_{od}$, determines the upper bound of the confidence interval, i.e. the upper bound of the approximated total surgery duration.
Combining equations \eqref{eq:mu-sigma}-\eqref{eq:overtimelog}, we obtain the following overtime constraints:
\begin{equation} \label{eq:overtimelognormcomplete}
    e^{\mathop{\mathbb{E}}[\tilde{k}_{od}] - \frac{1}{2}\ln \left( \frac{\text{Var}(\tilde{k}_{od})}{\mathop{\mathbb{E}}[\tilde{k}_{od}]^2}+1 \right) + \sqrt{\ln \left( \frac{\text{Var}(\tilde{k}_{od})}{\mathop{\mathbb{E}}[\tilde{k}_{od}]^2}+1 \right)}\Phi^{-1}(1-\alpha) } \leq C_{od},  \quad \forall  o \in \mathcal{O}, \quad d \in \mathcal{D}.
\end{equation}
We see that these overtime constraints are highly non-linear. Therefore, we introduce a prediction model to linearize these constraints.
This leads to the following overtime constraints, where $h_{od}(\mathop{\mathbb{E}}[\tilde{k}_{od}], \text{Var}(\tilde{k}_{od}))$ denotes the prediction model $h$ that predicts the $(1-\alpha)$-percentile of the approximated lognormal distribution of the total surgery duration in OR $o \in \mathcal{O}$ on day $d \in \mathcal{D}$.
The expected duration $\mathop{\mathbb{E}}[\tilde{k}_{od}]$ and the variance $\text{Var}(\tilde{k}_{od})$ of the total surgery duration on OR $o \in \mathcal{O}$ on day $d \in \mathcal{D}$ are the input variables of the model. 
\begin{equation}\label{eq:fnn_overtime}
    h_{od}(\mathop{\mathbb{E}}[\tilde{k}_{od}], \text{Var}(\tilde{k}_{od})) \leq C_{od}, \quad \forall o \in \mathcal{O}, \quad d \in \mathcal{D}.
\end{equation}

As mentioned in Section \ref{Sec:intro}, there are several regression models that are suitable to linearize constraints. 
\textcite{Fajemisin2024OptimizationSurvey} give an overview of the balance between explainability, complexity, required data and performance of predictive models used for constraint learning. It shows that, while neural networks lack in explainability and require many data points to train properly, they have the highest performance compared to regression trees and linear regression models. 
Given that the prediction model is used within hard constraints in our MILP model, performance is a critical aspect. Therefore, we select a neural network. Neural networks can be divided into three basic types of mapping networks \parencite{Meyer-Baese2014FoundationsNetworks} of which the feedforward neural network (FNN) approximates non-linear relationships based on a known input and output. 
Since we aim to approximate the Fenton-Wilkinson approximation using a known mean, variance and $(1-\alpha)$-percentile of the approximated lognormal distribution of the total surgery duration, we choose an FNN to linearize the constraints. 

Since the approximation of the $(1-\alpha)$-percentile is a regression problem, we apply the ReLU activation function as activation functions for the neurons of the FNN \parencite{Aggarwal2023NeuralLearning}.
To incorporate the FNN within an MILP model, each neuron needs to be expressed using a set of constraints. The ReLU activation function can be explicitly formulated using big-M constraints and binary auxiliary variables for each neuron \parencite{Maragno2025Mixed-IntegerLearning}. 
\textcite{Grimstad2019ReLUPrograms} provide an exact formulation of the ReLU activation function based on the work of \textcite{Fischetti2018DeepOptimization}. They present a general set of constraints for each layer to model a neural network consisting of ReLU activation functions.
An FNN with ReLU activation functions can also be implemented implicitly within an MILP model, provided that the used solver supports these type of constraints. Certain MILP solvers, such as Gurobi \parencite{GurobiOptimizationLLC2024GurobiLearning}, offer built-in functionality to embed trained machine learning models directly as constraints within a defined model. This allows us to incorporate the FNN within an MILP model without explicitly defining a set of constraints for each neuron. In this work, we apply the latter implicit approach to predict the $(1-\alpha)$-percentile of the approximated lognormal distribution of the total surgery duration.


\subsection{Constraint linearization with piecewise linear function} \label{sec:methodpiecewise}
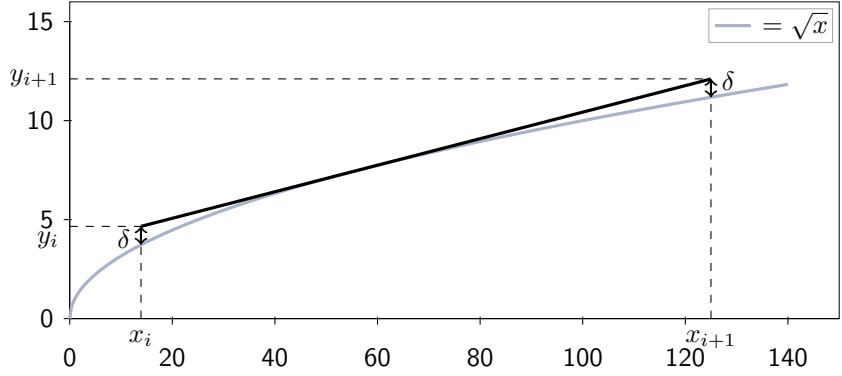
\begin{figure}[t]
\centering
\begin{tikzpicture}[
    scale=0.75,
    y=.35cm,
    x=.09cm,
    font=\sffamily\small
]

\draw (0,0) -- coordinate (x axis mid) (150,0);
\draw (0,0) -- coordinate (y axis mid) (0,15);
\draw[black, line width=0.3pt]
    (0,0) rectangle (150,16);
    
\foreach \x in {0,20,...,140}
    \draw (\x,1pt) -- (\x,-3pt)
        (\x,-10pt) node[anchor=north] {\x};

\foreach \y in {0,5,...,15}
    \draw (1pt,\y) -- (-3pt,\y)
        (-3pt,\y) node[anchor=east] {\y};

\draw[color=FNNcolor, very thick, domain=0:140, samples=400]
    plot (\x,{sqrt(\x)});

\draw[color=black, very thick, domain=13.89:125, samples=400]
    plot (\x,{3.729 + 0.0670*\x});

\draw[<->, thick]
    (13.89,3.73) -- (13.89,4.66)
    node[font=\small, pos=0.3, left] {$\delta$};

\draw[<->, thick]
    (125,11.18) -- (125,12.11)
    node[font=\small, pos=0.9, right] {$\delta$};

\draw[dashed] (0,4.66) -- (13.89,4.66);
\draw[dashed] (13.89,0) -- (13.89,3.73);
\draw[dashed] (0,12.11) -- (125,12.11);
\draw[dashed] (125,0) -- (125,11.18);

\draw (13.89,-0.1) node[font=\small, anchor=north] {$x_i$};
\draw (-0.1,4) node[font=\small, anchor=east] {$y_i$};

\draw (125,-0.1) node[font=\small, anchor=north] {$x_{i+1}$};
\draw (-0.1,12.11) node[font=\small, anchor=east] {$y_{i+1}$};

\begin{scope}[shift={(125.8,14.8)}]
    \draw[gray!70, fill=white, line width=0.3pt] 
        (-1.2,-0.85) rectangle (22.9,0.85);
    \draw[FNNcolor, very thick] (0,0) -- (7,0)
        node[right, text=black] {$=\sqrt{x}$};
\end{scope}

\end{tikzpicture}
\caption{Example of approximation of $\sqrt{x}$ on interval [13.5; 125.0] using a linear function.}
\label{fig:ex_breakpoints}
\end{figure}

\noindent
The approach to linearize the overtime constraints using piecewise linear functions is based on the work of \textcite{Schneider2020SchedulingResources}. For this method, we assume that the surgery durations are i.i.d., following a normal distribution, meaning that the total surgery duration also follows a normal distribution.  
We can then formulate overtime constraints \eqref{eq: generalovertime} as follows, where $\tilde{k}_{od}$ is normally distributed, i.e. $\tilde{k}_{od} \sim \mathcal{N}(\mu_{od}, \sigma^2_{od})$ for all ORs $o \in \mathcal{O}$ and days $d \in \mathcal{D}$: 
\begin{equation} \label{eq:normovertime}
    P(\tilde{k}_{od} \leq C_{od}) = \Phi \left( \frac{C_{od}-\mu_{od}}{\sigma_{od}} \right) \geq 1 - \alpha, \quad \forall  o \in \mathcal{O}, \quad d \in \mathcal{D}.
\end{equation}
These overtime constraints can be rewritten as:
\begin{equation} \label{eq:overtimenormaldist}
     \mu_{od} + \sigma_{od}\Phi^{-1}(1-\alpha) \leq C_{od},  \quad \forall  o \in \mathcal{O}, \quad d \in \mathcal{D},
\end{equation}
with $\mu_{od}$ and $\sigma^2_{od}$ defined as follows:
\begin{equation} \label{eq: mean_var_piecewise}
    \mu_{od} = \sum_{s \in \mathcal{S}} X_{sod}\mu_s \quad \text{and} \quad \sigma^2_{od} = \sum_{s \in \mathcal{S}} X_{sod}\sigma^2_s, \quad \forall  o \in \mathcal{O}, \quad d \in \mathcal{D}.
\end{equation}
Next, we combine equations \eqref{eq:overtimenormaldist} and \eqref{eq: mean_var_piecewise}, which leads to the following overtime constraints:
\begin{equation} \label{eq:overtimeconstrsqrt}
   \sum_{s \in \mathcal{S}} X_{sod}\mu_s + \Phi^{-1}(1-\alpha) \sqrt{\sum_{s \in \mathcal{S}} X_{sod}\sigma^2_s} \leq C_{od}, \quad \forall  o \in \mathcal{O}, \quad d \in \mathcal{D}.
\end{equation}
The square root function, which we denote as $f(x) = \sqrt{x}$, is a non-linear component within these constraints. 
To linearize this component, we apply a piecewise linear function with which we construct an approximation of the square root function, similar to the approach of \textcite{Schneider2020SchedulingResources}.
We approximate the function $f(x)$ on the interval $[x_{\text{min}}, x_{\text{max}}]$, where $x_{\text{min}}$ and $x_{\text{max}}$ denote the minimum and maximum possible value of $\sum_{s \in \mathcal{S}} X_{sod}\sigma^2_s$, respectively. 
Given that an OR might not have any surgeries scheduled, we set $x_{\text{min}}=0$. 
To determine $x_{\text{max}}$, we multiply the historical maximum number of surgeries scheduled in a single OR, denoted by $n$, by the highest observed surgery variance: 
\begin{equation}\label{eq:xmax-piecewise}
    x_{\text{max}} = n \cdot \max_{s \in \mathcal{S}}\{\sigma^2_{s}\}.
\end{equation}
Next, we split the interval $ \left[ 0, x_{\max} \right] $ into $b$ smaller subintervals, determined by breakpoints $i \in \mathcal{I}=\{0,1,\dots,b \} $, where $x_i$ denotes the value on the $x$-axis of breakpoint $i \in \mathcal{I}$.
We define $x_{0}=0$ and $x_b=x_{\max}$, and the other $x$-values mark the intersection points of the linear approximation functions.
As we do not want to underestimate the function $f(x)$, the approximation function must be greater than or equal to $f(x)$ for all $x \in \left[ 0, x_{\max} \right]$.
A higher number of breakpoints results in a more accurate approximation of the square root function. 
However, it also increases the number of variables within the model. 
Therefore, we introduce a maximum approximation error, denoted as $\Delta^{\text{max}}$. This $\Delta^{\text{max}}$ corresponds to the maximum allowable overestimated value of the approximation of the square root function at the breakpoints. 
Using this $\Delta^{\text{max}}$, we follow the method of \textcite{Schneider2020SchedulingResources} to determine the minimum number of breakpoints needed, their value $x_i$ for $i \in \{1,\dots,b-1\}$ and the corresponding function values $y_i$ for all breakpoints $x_i$ with $i \in \{0,1,\dots,b\}$.
The linear approximation function within each subinterval $\left[x_i, x_{i+1} \right]$ for $i = \{0,\dots,b-1\}$ is now given by the linear function connecting points $(x_i, y_i)$ and $(x_{i+1}, y_{i+1})$, which is by construction tangent to the square root function inside the considered interval.
This results in a certain error $\delta \leq \Delta^{\text{max}}$ that corresponds to the overestimation of the function value, which ultimately leads to $y_{i} = \sqrt{x_i} + \delta$.
Figure \ref{fig:ex_breakpoints} shows an example of an approximation of a section of $\sqrt{x}$ using a linear function.

Next, we apply the $\lambda$-formulation of \textcite{Bisschop2006AimmsModeling} to combine the piecewise linear functions. 
We assign a non-negative weight $\lambda_{odi} \in [0, 1]$ to each $i \in \mathcal{I}$, such that $\sum_{i \in \mathcal{I}} \lambda_{odi} = 1$ for each OR $o \in \mathcal{O}$ and each day $d \in \mathcal{D}$.
Then, we can rewrite overtime constraints \eqref{eq:overtimeconstrsqrt} as follows:
\begin{equation}\label{eq:plf_overtime}
   \sum_{s \in \mathcal{S}} X_{sod}\mu_s + \Phi^{-1}(1-\alpha) \sum_{i \in \mathcal{I}}\lambda_{odi}y_{i} \leq C_{od}, \quad \forall \text{ } o \in \mathcal{O}, \quad d \in \mathcal{D},
\end{equation}
\begin{equation}
    \sum_{i \in \mathcal{I}}\lambda_{odi}x_{i} = \sum_{s \in \mathcal{S}} X_{sod}\sigma^2_s, \quad \forall \text{ } o \in \mathcal{O}, \quad d \in \mathcal{D},
\end{equation}
\begin{equation}
    \sum_{i \in \mathcal{I}}\lambda_{odi} = 1, \quad \forall \text{ } o \in \mathcal{O}, \quad d \in \mathcal{D}.
\end{equation}
With this formulation, we have linearized the overtime constraints by approximating a normal distribution using piecewise linear functions.
The work of \textcite{Schneider2020SchedulingResources} provides a more detailed explanation of the methodology for determining the breakpoints and constructing the piecewise linear functions.

\subsection{Scenario-based modelling}\label{sec:scen-basedmodel}
\noindent
Lastly, we discuss scenario-based modelling as approach to incorporate the stochasticity of surgery durations.
In this approach, we create a finite set of scenarios that represent different possible combinations of surgery durations.
We define $\mathcal{L}$ as this set of scenarios with $|\mathcal{L}|$ the number of scenarios in set $\mathcal{L}$. 
Each scenario $l \in \mathcal{L}$ represents a joint realization of the stochastic surgery durations of all surgeries $s \in \mathcal{S}$.
To incorporate these stochastic durations, we assume that the surgery durations are lognormally distributed, i.e. $\tilde{w}_{s} \sim \text{LogNormal}(\mu_{s}, \sigma^2_{s})$ for all $s \in \mathcal{S}$. The scenarios $l \in \mathcal{L}$ are then generated by sampling realisations $w_{s}^{l}$ from these lognormal distributions. 

Note that the decision variable $X_{sod}$ is shared across all scenarios, and thus, is not scenario dependent.
To account for the overtime probability $\alpha$ in the general overtime constraints \eqref{eq: generalovertime}, we introduce chance constraints.
Let $y_{od}^{l} \in \{0, 1\}$ be a binary variable which is 1 if the total surgery duration in OR $o \in \mathcal{O}$ on day $d \in \mathcal{D}$ exceeds the total capacity $C_{od}$ under a certain scenario $l \in \mathcal{L}$, and 0 otherwise. 
When $y_{od}^{l}$ equals 1, the overtime constraint is relaxed with a large parameter, denoted as $M_{od}^{l}$. 
This leads to the following overtime constraints:
\begin{equation}\label{eq:constr_bigM}
    \sum_{s \in \mathcal{S}}X_{sod}w_s^{l} \leq C_{od} + M_{od}^{l}y_{od}^{l}, \quad \forall \text{ } l \in \mathcal{L}, \quad o\in \mathcal{O}, \quad d \in \mathcal{D}.    
\end{equation}
Together, constraints \eqref{eq:max_one}–\eqref{eq:overtime_samplemean} and constraints \eqref{eq:constr_bigM} define the complete constraint set for the scenario-based modelling approach.
Note that by introducing fixed parameter $w_s^l$ for the surgery duration in constraints \eqref{eq:constr_bigM}, we are dealing with an ILP model instead of an MILP model.

Constraints \eqref{eq:constr_bigM}, furthermore, show that parameter $M^{l}_{od}$ is dependant on scenario $l \in \mathcal{L}$, OR $o \in \mathcal{O}$ and day $d \in \mathcal{D}$, which ensures that the parameter is as tight as possible for each scenario-OR-day combination.
We compute the minimal feasible value of $M^{l}_{od}$ via auxiliary ILP models, which determine the maximum feasible total duration (including overtime) for each OR $o \in \mathcal{O}$, day $d \in \mathcal{D}$ and scenario $l \in \mathcal{L}$, denoted as $Q^{l}_{od}$.
This auxiliary ILP model uses the same generic constraints \eqref{eq:max_one} - \eqref{eq:overtime_samplemean} as the general problem statement, but has a different objective value function that determines $Q^l_{od}$:
\begin{equation}
   Q^{l}_{od} = \max \sum_{s \in \mathcal{S}}X_{sod}w_s^l, \quad \forall l \in \mathcal{L}, \quad o \in \mathcal{O}, \quad d \in \mathcal{D}. 
\end{equation}
With the maximum possible total duration $Q_{od}^l$, we determine $M_{od}^l$ as follows: 
\begin{equation}    M^{l}_{od} = Q^{l}_{od} - C_{od}, \quad \forall l \in \mathcal{L}, \quad o \in \mathcal{O}, \quad d \in \mathcal{D}. 
\end{equation}
Note that for a given day $d$ and scenario $l$, $M^{l}_{od}$ will be identical across any subset of ORs $o \in \mathcal{O}$ that have equal capacity $C_{od}$. In the current problem description, the ordering of ORs for a given day $d$ and scenario $l$ is irrelevant. We exploit this property to avoid running the ILP for every combination of $o \in \mathcal{O}, d \in \mathcal{D}$ and $l \in \mathcal{L}$. Therefore, we apply Algorithm \ref{alg:algorithmbasic} to determine $M^{l}_{od}$. 
Note that steps 2 and 3 are not combined, as it is important that the loop over set $\mathcal{O}$ is performed last.

\RestyleAlgo{ruled}
\begin{algorithm}
 \caption{$M_{od}^l$ algorithm}\label{alg:algorithmbasic}
\textbf{Input sets:} $\mathcal{S},\mathcal{O},\mathcal{D}, \mathcal{L}$. \\
\textbf{Input parameters:} $C_{od} \quad \forall o \in \mathcal{O}, \quad d \in \mathcal{D}$.\\
\textbf{Initiate:} $l_{\text{prev}}, C_{\text{prev}}, d_{\text{prev}}= \text{NaN}$ and $M_{od}^l= 0 \quad \forall l \in \mathcal{L}, \quad o \in \mathcal{O}, \quad d \in \mathcal{D}$. \\
\textbf{Algorithm steps:}
\begin{enumerate}
    \item Order ORs $o \in \mathcal{O}$ on capacity size $C_{od}$ for each $d \in \mathcal{D}$.
    
    \item \textbf{for} $l\in\mathcal{L}, d\in\mathcal{D}$:\\
    \item \textbf{for} $o\in\mathcal{O}$:\\
    \quad\textbf{if} $(l_{\text{prev}} = l)$ AND $(C_{\text{prev}} = C_{od})$ AND $(d_{\text{prev}} = d)$:\\
    \qquad \textbf{then} $M_{od}^{l} = M_{o-1,d}^{l}$ \\
    \quad\textbf{else:} solve the auxiliary ILP to determine $M_{od}^{l}$.\\
    \quad\textbf{Update:} $l_{\text{prev}} = l$,
    $C_{\text{prev}} = C_{od}$,
    $d_{\text{prev}} = d$.
\end{enumerate}
\KwResult{values $M_{od}^{l} \quad \forall \quad l \in \mathcal{L}, \quad o \in \mathcal{O}, \quad d \in \mathcal{D}$}
\end{algorithm}

To ensure that the predefined capacities $C_{od}$ are not exceeded for OR $o \in \mathcal{O}$ and day $d \in \mathcal{D}$ over the scenario space with a minimum probability of $(1-\alpha)$, we introduce 
\begin{equation}\label{eq:sbm_overtime}
    \sum_{l \in \mathcal{L}} y_{od}^{l} \leq \alpha|\mathcal{L}|, \quad  \forall o \in \mathcal{O}, \quad d \in \mathcal{D}.
\end{equation}
This enforces that for each OR $o \in \mathcal{O}$ on day $d \in \mathcal{D}$, the maximum number of scenarios that violate capacity $C_{od}$ under a certain decision $X_{sod}$ is at most $\alpha |\mathcal{L}|$, which is a predefined fraction of the total number of scenarios.
The creation of the scenarios $l \in \mathcal{L}$ and the determination of the number of scenarios $|\mathcal{L}|$ are explained in Section \ref{sec:createscenarios} and Section \ref{sec:resultscenarios}, respectively.

\section{Data description and processing}\label{sec:data_processing}
\noindent
To compare the three proposed methods, we use input data provided by the Erasmus MC. 
These data consist of 10,102 elective surgeries performed between July 2017 and December 2019.
Among these, 7,707 belong to the specialty cardiology (CAR), while 2,395 fall under the specialty otorhinolaryngology (ENT).
To use these data as input for the models, we first need to define model parameters and preprocess the data, which we describe in Section \ref{sec:processinputdata}. 
Furthermore, we need data to fit the FNN as part of the solution method that linearizes the overtime constraints using a neural network. 
We use a separate dataset to train this FNN to avoid bias in the performance. This separate dataset is provided by another Dutch academic medical centre, of which its accompanying data preprocessing steps are explained in Section \ref{sec:dataFNN}.
Lastly, we generate scenarios for the scenario-based modelling approach using the preprocessed data resulting from Section \ref{sec:processinputdata}. These steps are clarified in Section \ref{sec:createscenarios}.

\subsection{Input data models - model parameter settings \& preprocessing surgeries}\label{sec:processinputdata}
\noindent
In this section, we define model parameter settings and preprocess the dataset. 
First, we elaborate on the general parameter settings in Section \ref{sec:generalparameters}, which are based on contextual hospital insights.
Next, we perform general preprocessing steps in Section \ref{sec:generalprocess}.
These general steps are required to define the approach-specific parameters. We describe the parameters in Section \ref{sec:specificparameters}, after which we generate four different input datasets in Section \ref{sec:createinputdatasets}. 
These four datasets make it possible to compare alternative scheduling situations across our different approaches.

\subsubsection{General parameter settings}\label{sec:generalparameters}
\noindent
First, we assume that elective surgeries can only be scheduled on weekdays. 
Therefore, in our research, each week consists of 5 days. 
In general, elective surgery schedules are drawn up and evaluated on a weekly basis.
Thus, we set our scheduling horizon to one week. 
Next, we define the numbers of available ORs for both cardiology and ENT, which are based on the hospital insights from 2019 and visualized in Figure \ref{fig:avail_ORs}. 
For cardiology, 4 ORs are available on the first and third day of the week, while 5 ORs are available on the other days. These ORs are all open from 08:00 to 16:30.
For ENT, each day, 2 ORs are available from 08:00 to 16:30. There is one exception on the fifth day of the week, when one OR is opened until 20:00.
Lastly, we set the overtime probability $\alpha$ equal to 0.15, which leads to a threshold of 0.85 in overtime constraints \eqref{eq: generalovertime}. This now states that the probability of the total surgery duration $\tilde{k}_{od}$ being less than or equal to the maximum OR capacity $C_{od}$ must be at least 0.85. This condition holds for each OR $o \in \mathcal{O}$ and each day $d \in \mathcal{D}$.

\begin{figure}[t]
    \centering
    \begin{tikzpicture}
\matrix (m) [matrix of nodes,
    font=\footnotesize, 
    nodes={draw, align=center, outer sep=0pt,
           anchor=center}, 
    column sep=-\pgflinewidth,
    row sep=0pt 
]{
      & |[draw=none]| \textbf{Day 1} & |[draw=none]| \textbf{Day 2} & |[draw=none]|\textbf{Day 3} & |[draw=none]|\textbf{Day 4} & |[draw=none]|\textbf{Day 5} \\
     |[draw=none, minimum width = 1.5cm]| OR 1 
     & \node[fill=CARcolor, minimum width=2.4cm, minimum height=0.8cm]{CAR\\08:00 - 16:30};
     & \node[fill=CARcolor, minimum width=2.4cm, minimum height=0.8cm]{CAR\\08:00 - 16:30};
     & \node[fill=CARcolor, minimum width=2.4cm, minimum height=0.8cm]{CAR\\08:00 - 16:30};
     & \node[fill=CARcolor, minimum width=2.4cm, minimum height=0.8cm]{CAR\\08:00 - 16:30};
     & \node[fill=CARcolor, minimum width=2.4cm, minimum height=0.8cm]{CAR\\08:00 - 16:30};\\
    |[draw=none, minimum width = 1.5cm]| OR 2  
     & \node[fill=CARcolor, minimum width=2.4cm, minimum height=0.8cm]{CAR\\08:00 - 16:30};
     & \node[fill=CARcolor, minimum width=2.4cm, minimum height=0.8cm]{CAR\\08:00 - 16:30};
     & \node[fill=CARcolor, minimum width=2.4cm, minimum height=0.8cm]{CAR\\08:00 - 16:30};
     & \node[fill=CARcolor, minimum width=2.4cm, minimum height=0.8cm]{CAR\\08:00 - 16:30};
     & \node[fill=CARcolor, minimum width=2.4cm, minimum height=0.8cm]{CAR\\08:00 - 16:30};\\
    |[draw=none, minimum width = 1.5cm]| OR 3  
     & \node[fill=CARcolor, minimum width=2.4cm, minimum height=0.8cm]{CAR\\08:00 - 16:30};
     & \node[fill=CARcolor, minimum width=2.4cm, minimum height=0.8cm]{CAR\\08:00 - 16:30};
     & \node[fill=CARcolor, minimum width=2.4cm, minimum height=0.8cm]{CAR\\08:00 - 16:30};
     & \node[fill=CARcolor, minimum width=2.4cm, minimum height=0.8cm]{CAR\\08:00 - 16:30};
     & \node[fill=CARcolor, minimum width=2.4cm, minimum height=0.8cm]{CAR\\08:00 - 16:30};\\
    |[draw=none, minimum width = 1.5cm]| OR 4  
     & \node[fill=CARcolor, minimum width=2.4cm, minimum height=0.8cm]{CAR\\08:00 - 16:30};
     & \node[fill=CARcolor, minimum width=2.4cm, minimum height=0.8cm]{CAR\\08:00 - 16:30};
     & \node[fill=CARcolor, minimum width=2.4cm, minimum height=0.8cm]{CAR\\08:00 - 16:30};
     & \node[fill=CARcolor, minimum width=2.4cm, minimum height=0.8cm]{CAR\\08:00 - 16:30};
     & \node[fill=CARcolor, minimum width=2.4cm, minimum height=0.8cm]{CAR\\08:00 - 16:30};\\
    |[draw=none, minimum width = 1.5cm]| OR 5  
     & 
     & \node[fill=CARcolor, minimum width=2.4cm, minimum height=0.8cm]{CAR\\08:00 - 16:30};
     & 
     & \node[fill=CARcolor, minimum width=2.4cm, minimum height=0.8cm]{CAR\\08:00 - 16:30};
     & \node[fill=CARcolor, minimum width=2.4cm, minimum height=0.8cm]{CAR\\08:00 - 16:30};\\
     |[draw=none, minimum height=0pt]| \vspace{0.002cm} 
     & |[draw=none]| \vspace{0.002cm} 
     & |[draw=none]| \vspace{0.002cm} 
     & |[draw=none]| \vspace{0.002cm} 
     & |[draw=none]| \vspace{0.002cm} 
     & |[draw=none]| \vspace{0.002cm} \\
     |[draw=none, minimum width = 1.5cm]| OR 6  
     & \node[fill=ENTcolor, minimum width=2.4cm, minimum height=0.8cm]{ENT\\08:00 - 16:30};
     & \node[fill=ENTcolor, minimum width=2.4cm, minimum height=0.8cm]{ENT\\08:00 - 16:30};
     & \node[fill=ENTcolor, minimum width=2.4cm, minimum height=0.8cm]{ENT\\08:00 - 16:30};
     & \node[fill=ENTcolor, minimum width=2.4cm, minimum height=0.8cm]{ENT\\08:00 - 16:30};
     & \node[fill=ENTcolor, minimum width=2.4cm, minimum height=0.8cm]{ENT\\08:00 - 16:30};\\
     |[draw=none, minimum width = 1.5cm]| OR 7  
     & \node[fill=ENTcolor, minimum width=2.4cm, minimum height=0.8cm]{ENT\\08:00 - 16:30};
     & \node[fill=ENTcolor, minimum width=2.4cm, minimum height=0.8cm]{ENT\\08:00 - 16:30};
     & \node[fill=ENTcolor, minimum width=2.4cm, minimum height=0.8cm]{ENT\\08:00 - 16:30};
     & \node[fill=ENTcolor, minimum width=2.4cm, minimum height=0.8cm]{ENT\\08:00 - 16:30};
     & \node[fill=ENTcolordark, minimum width=2.4cm, minimum height=0.8cm]{ENT\\08:00 - 20:00};\\
};
\end{tikzpicture}
    \caption{Availability of ORs for both cardiology (CAR) and ENT. The colour intensity indicates different time slots, also mentioned below the specialty.}
    \label{fig:avail_ORs}
\end{figure}

\subsubsection{General data preprocessing} \label{sec:generalprocess}
\noindent
We start with the previously mentioned dataset comprising 10,102 elective surgeries: 7,707 cardiology and 2,395 ENT procedures.
First, we exclude surgeries containing registration errors by removing all surgeries with a negative surgery duration or a duration of 0 minutes. 
Subsequently, we remove any surgery with a duration longer than 12 hours, which is the maximum plannable duration for elective surgeries and matches the defined model in Section \ref{Sec:problemdescr}.
This leaves us with 7,651 cardiology procedures consisting of 143 different surgery types and 2,323 ENT procedures with 112 different surgery types.
Each surgery type corresponds to a distinct category of surgical procedure.
The set of surgery types is defined by set $\mathcal{T}$ and the subset of surgeries having surgery type $t \in \mathcal{T}$ as main surgical procedure is given by $\mathcal{S}_t \subseteq \mathcal{S}$.
Next, we determine the sample mean $\bar{w}_t$ for each surgery type $t \in \mathcal{T}$ based on surgeries $s \in \mathcal{S}_t$. Subsequently, we determine parameters $\mu_t^{(N)}$ and $\sigma_t^{(N)}$, belonging to the normal distribution describing the surgery duration distribution of surgery type $t \in \mathcal{T}$, and $\mu_t^{(LN)}$ and $\sigma_t^{(LN)}$, belonging to the lognormal distribution describing the surgery duration distribution of surgery type $t \in \mathcal{T}$. If $|\mathcal{S}_t|=1$ for $t \in \mathcal{T}$, we set $\sigma_t$ to 0 for both distributions.
Lastly, we project the parameters of each surgery type $t \in \mathcal{T}$ to all surgeries $s \in \mathcal{S}_t$. 
Hence, for all surgeries $s \in \mathcal{S}_t$ belonging to surgery type $t \in \mathcal{T}$, we set $\bar{w}_s$, $\mu_s^{(N)}$, $\sigma_s^{(N)}$, $\mu_s^{(LN)}$, and $\sigma_s^{(LN)}$  equal to $\bar{w}_t$, $\mu_t^{(N)}$, $\sigma_t^{(N)}$, $\mu_t^{(LN)}$, and $\sigma_t^{(LN)}$, respectively.
Note that sample mean $\bar{w}_s$ is equal to parameter $\mu_s^{(N)}$ for each $s \in \mathcal{S}$.

\subsubsection{Approach specific parameter settings} \label{sec:specificparameters}
\noindent
Next, we define approach specific parameters. 
Only the constraint linearization with piecewise linear functions approach contains specific parameters.
For this approach, we have to determine the interval $[x_{\text{min}}, x_{\text{max}}]$.
We have already stated in Section \ref{sec:methodpiecewise} that $x_{\text{min}}$ equals 0, since a schedule can contain an empty OR. 
We define $x_{\text{max}}$ using equation \eqref{eq:xmax-piecewise}. Based on the preprocessed dataset, the maximum number of surgeries scheduled in a single OR equals 8 for both specialties and the highest observed surgery duration variance over both specialties equals 54,035. This leads to an $x_{\text{max}}$ of 432,280.
For this approach, we also define the maximum approximation error $\Delta^{\text{max}}$, which we set equal to 1 minute.
This results in a total of 19 breakpoints to approximate the square root function in constraints \eqref{eq:overtimeconstrsqrt} with $\delta = 0.965$.

\subsubsection{Generate separate input datasets} \label{sec:createinputdatasets}
\noindent
Finally, we construct four datasets, each with a planning horizon of five days. Two datasets contain cardiology surgeries, and two include ENT surgeries.
Since the number of available ORs is based on data from 2019, we use surgeries with a due date in 2019.
To create a dataset, we pick an arbitrary week in 2019 and select all surgeries that are on the waiting list or appear on the waiting list that specific week. 
Next, we determine release date $r_s$ and due date $q_s$ for all surgeries $s \in \mathcal{S}$.
If the surgery is on the waiting list on the first day of the week, we set $r_s$ equal to the first day of the planning horizon, which is 0. 
For the remaining surgeries, the release date $r_s$ corresponds to the next consecutive day they appear on the list. 
If due date $q_s$ of surgery $s \in \mathcal{S}$ is missing, we assume the surgery has no due date. In that case, we set $q_s$ equal to the maximum due date among all surgeries within the dataset plus 1.
This leads to the four final datasets with their specifics depicted in Table \ref{tab:datasets-settings}. 
In this table, we also show the mean standard deviation of the surgeries, which captures the average variability of surgeries within a dataset and is quantified as 
$ \frac{1}{|\mathcal{S}|}\sum_{s\in\mathcal{S}}\sigma_s^{(N)}$.   
It shows that there is on average a higher variability in surgery duration within the ENT datasets compared to both cardiology datasets.

\begin{table}[t]
    \centering
    \caption{Overview 4 different datasets used as model input}
    \label{tab:datasets-settings} \begin{tabular}{lcccc}
    \toprule
    & \multicolumn{2}{c}{\textbf{Cardiology}} & \multicolumn{2}{c}{\textbf{ENT}}\\
    \cmidrule(lr){2-3}
    \cmidrule(lr){4-5}
    
      & Dataset 1 & Dataset 2 & Dataset 1 & Dataset 2 \\
    \midrule
      \# Surgeries: $|\mathcal{S}|$ & 216 & 158 & 137 & 52 \\
        \# Different surgery types: $|\mathcal{T}|$ & 50 & 40 & 45 & 31 \\
    \# Surgeries with due date in horizon: $p_s=1$ & 13 & 12 & 5 & 4 \\
      
    \# Surgeries with release date: & & & & \\
    \hspace{1em} $r_s = 0$ & 184 & 136 & 129 & 48 \\
    \hspace{1em} $r_s = 1$ & 18 & 6 & 5 & 1 \\
    \hspace{1em} $r_s = 2$ & 3 & 6 & 0 & 2\\
    \hspace{1em} $r_s = 3$ & 10 & 6 & 2 & 1 \\
    \hspace{1em} $r_s = 4$ & 1 & 4 & 1 & 0\\

  Total sample mean duration (min): $\sum_{s \in \mathcal{S}}\bar{w}_s$ & 24,191 & 17,543 & 18,808 & 6,279  \\
    Total OR capacity (min) & 11,730 & 11,730 & 5,310 & 5,310 \\ 
    Mean SD (min) & 45.4 & 46.1 & 62.9 & 54.2 \\
    \bottomrule
    \end{tabular}
\end{table}

Lastly, we determine for each surgery type that contains at least 5 surgeries in the dataset whether the lognormal or normal distribution is a better fit for surgery duration.
We apply the Akaike Information Criterion (AIC), which is a measure used to compare statistical models by balancing goodness of fit with model complexity. Lower AIC values indicate a better trade-off between accuracy and simplicity.
Table \ref{tab:dist_surgery_types} presents the number of surgery types that are better fitted by either a normal or lognormal distribution for each dataset.
In general, the lognormal distribution is a better fit for the surgery types.
Furthermore, the cardiology datasets contain more surgery types for which the normal distribution is a better fit compared to the ENT datasets.

\begin{table}[!t]
    \centering
     \caption{Number of surgery types in the dataset that better fit a lognormal distribution versus a normal distribution.}
    \begin{tabular}{lcc}
    \toprule
     Dataset  & Lognormal distribution  & Normal distribution  \\
     \midrule
     Cardiology 1    & 28 & 17 \\
     Cardiology 2 & 22 & 13\\
     ENT 1 & 30 & 9\\
     ENT 2 & 23 & 5\\
     \bottomrule
    \end{tabular}
   \label{tab:dist_surgery_types}
\end{table}

\subsection{Data preprocessing to fit FNN} \label{sec:dataFNN}
\noindent
To capture the Fenton-Wilkinson approximation with an FNN, we need surgical data where the mean, variance and $85^{\text{th}}$ percentile of the approximated total surgery duration are known.
As previously mentioned, we use a separate dataset to fit the FNN to verify if the model is generalizable and to avoid bias.
These data belong to a different academic medical centre in the Netherlands and consist of 7,767 surgeries performed between January 2015 and March 2016. We preprocess and structure these data through the following steps.

\subsubsection{Determine expectation and variance using lognormal parameters}
\noindent
Our first goal is to determine the expectation and variance of the surgeries based on their surgery type. We start with grouping the surgeries on their surgery type, which results in 1,535 different surgery types. The smallest groups contain one surgery and the largest group contains 177 surgeries.
Next, we estimate the parameters of the lognormal distribution based on each surgery type.
In theory, the parameters of the lognormal distribution can be estimated with only 2 data points. However, this does not necessarily result in a reliable estimation. 
To obtain a more accurate estimation, we only include surgery types containing at least 30 surgeries. 
After filtering, 35 surgery types remain, comprising a total of 2,074 surgeries.
For each of these surgeries, we determine the parameters $\mu_s^{(LN)}$ and $\sigma_s^{(LN)}$ belonging to the lognormal distribution of their corresponding surgery type.
Using these parameters and equation \eqref{eq:exp_var_lognormal}, we determine the expectation $\mathop{\mathbb{E}}[\tilde{w}_s]$ and variance $\text{Var}(\tilde{w}_s)$ for each surgery $s \in S$. This means that each surgery with the same surgery type has the same expectation and variance. With these variables, we can resample the data and create new data points in the next step. 

\subsubsection{Resample data}
\noindent
To create a sufficient number of data points to fit an FNN, we resample the data. 
One important aspect is to create a prediction region larger than the region restricted by the overtime constraints to ensure that the FNN creates reliable predictions. 
Therefore, we create new data points by generating all possible surgery type combinations with replacement of sizes 1 through 6. This ensures that there will be data points with an approximated total surgery duration far outside the maximum capacity of the ORs.
For each of these combinations, we calculate the combined expectation and variance using equations \eqref{eq:totmean} and \eqref{eq:totvar}. 
This results in a total of 4,496,387 data points with an expected total surgery duration ranging from 42 minutes to 1,711 minutes.

\subsubsection{Finetuning the dataset and determine 85th percentile}
\noindent
To increase model performance, we remove data points that have outliers in either the feature `expectation' or `variance'. 
We consider an outlier to be a data point where the value for one of the features is less than $\mu - 3\sigma$ or greater than $\mu + 3\sigma$, where $\mu$ is that feature's mean and $\sigma$ its standard deviation.
This leads to a dataset containing 4,448,511 data points.   
Note that in the previous step, the range of the total expected surgery duration does not start at 0 minutes. However, it is possible to have an OR without any surgeries scheduled on it. To make sure that the FNN also fits these zero-valued data points properly, we add additional data points, each with a mean and variance of 0. These data points make up $1\%$ of the total dataset size, leading to a dataset containing 4,492,996 data points.

Lastly, we determine the $85^{\text{th}}$ percentile of the distribution for each data point using equation \eqref{eq:overtimelognormcomplete} with an $\alpha$ of 0.15. 
These steps result in the final dataset that we use to fit an FNN, where the expectation and variance of each data point are used as input features and the $85^{\text{th}}$ percentile of the approximated lognormal distribution of the total surgery duration is used as output feature. 
We elaborate on the FNN design in Section \ref{Sec:FNNmodeldesign}.

\subsection{Create scenarios for scenario-based modelling} \label{sec:createscenarios}
\noindent
The scenario-based modelling approach requires scenario sets for each separate input dataset. 
To construct these scenarios, we assume that the surgery durations are lognormally distributed. 
Therefore, for each $s \in \mathcal{S}$, we use parameters $\mu_s^{(LN)}$ and $\sigma_s^{(LN)}$, which we determined in Section \ref{sec:generalprocess}. 
We generate a set of scenarios $\mathcal{L}$, where each scenario $l \in \mathcal{L}$ corresponds to a vector containing sampled surgery durations $w_s^l$ for all $s \in \mathcal{S}$. These samples are generated by drawing values from the lognormal distribution with parameters $\mu_s^{(LN)}$ and $\sigma_s^{(LN)}$. 
Next, we determine the number of scenarios needed, defined as $|\mathcal{L}|$. 
Generating thousands of scenarios is a relatively efficient task. However, solving ILP models with thousands of constraints is computationally expensive. 
Therefore, choosing the smallest number of scenarios that is representative for the uncertainty of surgery durations is essential.
To determine the most suitable number of scenarios, we compare different sizes for $|
\mathcal{L}|$ as follows.

First, we generate 2,000 scenarios for all four datasets based on the lognormal distribution of the surgeries as mentioned above.
Subsequently, we use a scenario reduction technique called $k$-medoids clustering \parencite{Jin2011K-MedoidsClustering} to reduce the number of scenarios in each set while preserving variability. This technique creates a predefined number of clusters of similar scenarios and represents each cluster with a single scenario, called a medoid. 
These medoids become the reduced scenario set and allow a balance between accuracy and computational efficiency.
We use this technique to create sets containing scenarios with sizes 10 to 250 with regular increments of 10. 
To determine which of these sizes is the most suitable, we test the performance of the ILP models with these different number of scenarios and choose the smallest number of scenarios for which the outcome of the model is stable. 
We perform this analysis on all four input datasets.
The final number of scenarios is chosen in Section \ref{sec:resultscenarios}. 

\section{Experimental results}\label{Sec:results}
\noindent
In this section, we evaluate the three different approaches and present the results.
First, we design the FNN and determine the number of scenarios needed for scenario-based modelling in Sections \ref{Sec:FNNmodeldesign} and \ref{sec:resultscenarios}, respectively. 
We assess the performance of the approaches by comparing their generated schedules in Section \ref{sec:comp_sol_methods}, running a simulation study in Section \ref{sec:sim_study}, and evaluating schedule feasibility across all approaches in Section \ref{sec:feas_check}.
Lastly, we use real-life overtime cases as input for each approach to evaluate how they handle these cases in Section \ref{sec:real-life-cases}.
We solve the models using Python version 3.12 and Gurobi version 12.0.1. All experiments are performed on a PC with an Intel Core i7-13700H 2.40 gigahertz with 32 gigabytes RAM.

\subsection{FNN model design} \label{Sec:FNNmodeldesign}
\noindent
To linearize overtime constraints \eqref{eq:overtimelognormcomplete}, we make use of an FNN for the approximation of the total surgery duration. As previously stated, this approximation of the $85^{\text{th}}$ percentile of the total surgery duration probability function is a regression problem. Therefore, we use the ReLU activation function for the hidden layers and the mean squared error as loss function \parencite{Aggarwal2023NeuralLearning}. Since we are dealing with a duration prediction, which is non-negative, the activation for the output layer is a ReLU activation function as well. 
The input layer consists of two neurons (the approximated mean and variance of the total surgery duration), while the output layer consists of one neuron (the prediction of the $85^{\text{th}}$ percentile). 
We determine two parameters related to the size of the network, namely the number of hidden layers and the number of neurons, while also fine-tuning two parameters related to the training process of the FNN, which are the learning rate and batch size.

\textcite{Uzair2020EffectsNetworks} state that choosing the number of hidden layers is a balance between accuracy and computational efficiency. They also state that, while choosing the number of hidden layers and neurons is problem specific, excessive increase of these components can lead to overfitting. Furthermore, the number of binary variables needed in MILP models is proportional to the size of the neural network and smaller networks reduce the computational burden of solving the MILP model \parencite{Fajemisin2024OptimizationSurvey}. Since we are dealing with a relatively small problem (two input neurons and one output neuron), but are trying to fit a  non-linear relationship, we train the model with 1 to 3 hidden layers with each hidden layer containing one of the following numbers of neurons: 2, 4, 6, or 8. With these settings, we create a balance between accuracy and computational efficiency, while avoiding overfitting.
For the learning rate, commonly used values range from 0.1 to 0.001, normally chosen in logarithmic steps to balance convergence speed and stability \parencite{Skansi2018FeedforwardNetworks}. Therefore, we fit the model using a learning rate of 0.1, 0.01 and 0.001.
The batch size refers to the number of training samples used in one pass during training. A smaller batch size leads to more frequent updates, while a large batch size leads to faster computation with the risk of overfitting. For the batch size, we use the following sizes: 32, 64 and 128.
Lastly, we split the dataset into a train, validation, and test set with a split ratio of 0.7, 0.15, 0.15, respectively.

After fitting, the best performing model consists of 3 hidden layers with 8 neurons in each layer, and is trained with a learning rate of 0.01 and a batch size of 64. However, the model with 2 hidden layers and 8 neurons (learning rate: 0.01, batch size: 32) and the model with 3 hidden layers with 6 neurons in each layer (learning rate: 0.01, batch size: 64) have very similar performances. The performance measurements for each of these three models are visible in Table \ref{tab:performance_FNN}. Note that the mean absolute error does not change for the train, validation and test set for each model. This is due to very small differences that are not visible as a result of rounding.
All other possible parameter combinations resulted in a mean absolute error greater than 1 minute for the training, validation, and test sets. Since the FNN is used as part of a hard constraint within the models, accuracy is critical. As a result, we do not consider these configurations.
Each neuron in an FNN leads to an additional set of constraints within MILP models, which makes smaller models more desirable to incorporate. 
Therefore, we choose to use the FNN model with 2 layers and 8 neurons for further computations, which has the smallest total number of neurons.

\begin{table}[t]
    \centering
    \small
    \begin{threeparttable}
    \caption{Performance metrics of the three best performing 
    FNNs.} \label{tab:performance_FNN}
\begin{tabular}{l p{1.8cm} p{2.2cm} p{1.8cm} p{2.2cm} p{1.8cm} p{2.2cm}}
        \toprule
        &  \multicolumn{2}{c}{\textbf{3 layers, 8 neurons}} &   \multicolumn{2}{c}{\textbf{2 layers, 8 neurons}} & \multicolumn{2}{c}{\textbf{3 layers, 6 neurons}} \\
        \cmidrule(lr){2-3} \cmidrule(lr){4-5} \cmidrule(lr){6-7}
        \textbf{Set}  & MAE (min) & MaxAE (min) & MAE (min) & MaxAE (min) & MAE (min) & MaxAE (min) \\ 
        \midrule 
        Train set &  0.14 & 4.90 & 0.22 & 10.16 & 0.29 & 8.33 \\
        Validation set &  0.14 & 5.24 & 0.22 & 6.71 & 0.29 & 7.23 \\
         Test set &  0.14 & 4.82 & 0.22 & 8.04 & 0.29 & 7.10 \\
        \bottomrule 
    \end{tabular} 
    \begin{tablenotes}
\footnotesize
\item[] MAE: mean absolute error; MaxAE: max absolute error
\end{tablenotes}
\end{threeparttable}
\end{table}

\subsection{Number of scenarios scenario-based modelling}\label{sec:resultscenarios}
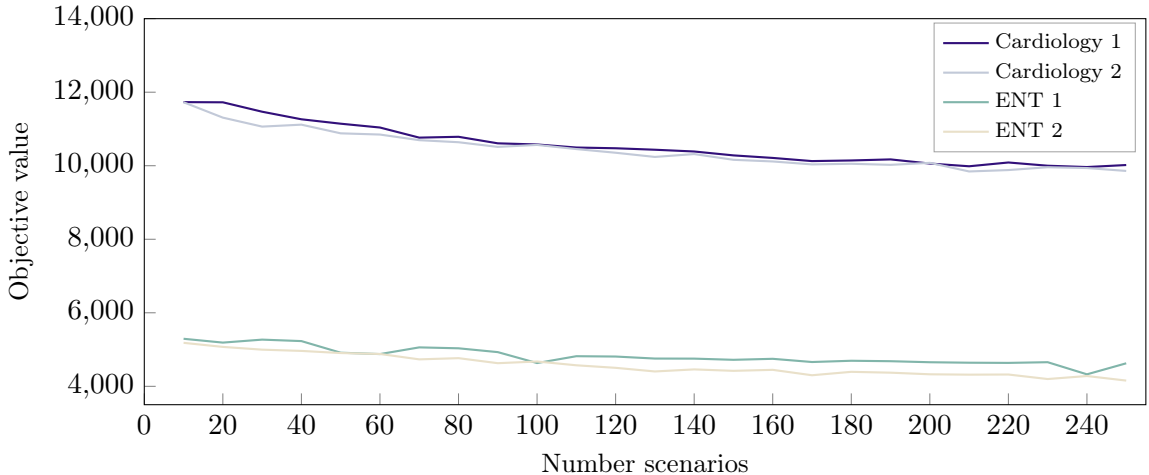
\begin{figure}[t]
\centering
\begin{tikzpicture}
\begin{axis}[
    width=0.9\textwidth,
    height=6.7cm,
    xlabel={\small Number scenarios},
    ylabel={\small Objective value},
    ymin=3500, ymax=14000,
    xmin = 0, xmax = 255,
    legend pos=north east,
        legend style={
        font=\scriptsize,
        legend cell align=left,
        at={(0.99,0.99)},
        anchor=north east,
        draw=gray!70
    },   
    scaled ticks=false,
    xtick pos=bottom, 
    ytick pos=left    
]

\pgfplotstableread[row sep=\\]{
scen car1 car2 ent1 ent2 \\
10	11730	11727	5293	5183 \\
20	11725	11308	5189	5072\\
30	11469	11064	5270	4997\\
40	11263	11117	5230	4963\\
50	11143	10881	4914	4905\\
60	11039	10848	4877	4880\\
70	10764	10695	5059	4733\\
80	10787	10639	5034	4766\\
90	10609	10512	4930	4629\\
100	10578	10566	4633	4675\\
110	10495	10453	4819	4572\\
120	10475	10353	4810	4503\\
130	10435	10239	4755	4404\\
140	10387	10317	4753	4460\\
150	10281	10161	4721	4421\\
160	10212	10116	4748	4448\\
170	10126	10036	4661	4301\\
180	10143	10053	4695	4394\\
190	10172	10024	4683	4371\\
200	10059	10078	4655	4327\\
210	9985	9844	4642	4317\\
220	10088	9882	4637	4321\\
230	10000	9959	4656	4198\\
240	9964	9937	4326	4280\\
250	10019	9858	4626	4156\\
}\mydata

\addplot[thick, CAR1color]   table[x=scen,y=car1]{\mydata};
\addplot[thick, CARcolor]  table[x=scen,y=car2]{\mydata};
\addplot[thick, ENT1color] table[x=scen,y=ent1]{\mydata};
\addplot[thick, ENTcolor] table[x=scen,y=ent2]{\mydata};

\legend{Cardiology 1, Cardiology 2, ENT 1, ENT 2}

\end{axis}
\end{tikzpicture}
\caption{Objective function value for different number of scenarios presented for each dataset.}
\label{fig:nr_scenarios}
\end{figure}

\noindent
To determine an appropriate number of scenarios, we evaluate the model's performance presented in Section \ref{sec:scen-basedmodel} using scenario set sizes $|\mathcal{L}|$ ranging from 10 to 250 with an increment of 10. For each scenario size, we execute the scenario-based model with a time limit of 12 hours.
Figure \ref{fig:nr_scenarios} shows the objective value for all four datasets under these different scenario sizes. 
As shown, the objective value stabilizes at approximately 210 scenarios for both cardiology datasets and 170 scenarios for both ENT datasets. These differences in number of scenarios arise from the varying input dimensions between the datasets. Consequently, for further computations, we use 210 and 170 scenarios for the cardiology and ENT settings, respectively.

\begin{table}[!htbp]
\caption{Performance of three approaches with different computation times}
\label{tab:results_models}
\centering

\begin{subtable}{\textwidth}
    \centering
    \small
    \caption{Performance Cardiology dataset 1} \label{tab:results_models_CAR1}
\begin{tabular}{p{3.5cm} p{1.3cm}p{1.3cm} p{1.3cm} p{1.3cm}p{1.3cm}p{1.3cm}}
        \toprule
        &  \multicolumn{3}{c}{\textbf{5 minutes}} &   \multicolumn{3}{c}{\textbf{1 hour}} \\
        \cmidrule(lr){2-4} \cmidrule(lr){5-7} 
        \textbf{KPI}  & FNN & PLF & SBM & FNN & PLF & SBM  \\ 
        \midrule 
        Objective value &  9,912 & 9,901& 9,408 &9,944 & 9,923 & 9,493\\
        \hangindent=1em \hangafter=0 Total duration: $\sum_{s \in \mathcal{S}}\bar{w}_s$  &  9,906 & 9,893& 9,399 &9,937 & 9,916 & 9,484\\
        \hspace{1em}Priority: $\sum_{s \in \mathcal{S}} \frac{1}{q_s+1}$ &  6.9 & 7.5& 8.9 &7.0 & 7.0 & 8.8 \\
        OR utilization \textsuperscript{1} (\%) & 84.5 & 84.3& 80.1 &84.7 & 84.5 & 80.9\\
        \# Surgeries &  91 & 97& 122 &92 & 92 & 122 \\
        Optimality gap (\%) &  1.05 & 4.78& 24.78 &0.62 & 2.01 & 23.66 \\
        \bottomrule 
    \end{tabular} 
\end{subtable}

\vspace{1em} 

\begin{subtable}{\textwidth}
    \centering
    \small
    \caption{Performance Cardiology dataset 2} \label{tab:results_models_CAR2}
\begin{tabular}{p{3.5cm} p{1.3cm}p{1.3cm} p{1.3cm} p{1.3cm}p{1.3cm}p{1.3cm}}
        \toprule
        &  \multicolumn{3}{c}{\textbf{5 minutes}} &   \multicolumn{3}{c}{\textbf{1 hour}} \\
        \cmidrule(lr){2-4} \cmidrule(lr){5-7} 
        \textbf{KPI}  & FNN & PLF & SBM & FNN & PLF & SBM  \\ 
        \midrule 
        Objective value &  9,774 & 9,654& 9,317 &9,801 & 9,654 & 9,408\\
        \hangindent=1em \hangafter=0 Total duration: $\sum_{s \in \mathcal{S}}\bar{w}_s$ &  9,765 & 9,646 & 9,308 &9,791 & 9,646 & 9,398\\
        \hspace{1em}Priority: $\sum_{s \in \mathcal{S}} \frac{1}{q_s+1}$ &  9.0 & 8.4& 9.8 &9.0 & 8.4 & 9.7\\
        OR utilization \textsuperscript{1} (\%) & 83.2 & 82.2& 79.4 &83.5 & 82.2 & 80.2 \\
        \# Surgeries & 96 & 95& 111 &96 & 95 & 112 \\
        Optimality gap (\%) &  0.99 & 3.27& 26.00 &0.70 & 2.89 & 24.79 \\
        \bottomrule 
    \end{tabular} 
\end{subtable}

\vspace{1em} 

\begin{subtable}{\textwidth}
    \centering
    \small
    \caption{Performance ENT dataset 1} \label{tab:results_models_ENT1}
\begin{tabular}{p{3.5cm} p{1.3cm}p{1.3cm} p{1.3cm} p{1.3cm}p{1.3cm}p{1.3cm}}
        \toprule
        &  \multicolumn{3}{c}{\textbf{5 minutes}} &   \multicolumn{3}{c}{\textbf{1 hour}} \\
        \cmidrule(lr){2-4} \cmidrule(lr){5-7} 
        \textbf{KPI}  & FNN & PLF & SBM & FNN & PLF & SBM  \\ 
        \midrule 
        Objective value &  4,544 & 4,528& 4,356 & 4,551 & 4,539 & 4,620\\
        \hangindent=1em \hangafter=0 Total duration: $\sum_{s \in \mathcal{S}}\bar{w}_s$ &  4,540 & 4,525& 4,351 &4,547 & 4,531 & 4,616\\
        \hspace{1em}Priority: $\sum_{s \in \mathcal{S}} \frac{1}{q_s+1}$ &  3.6 & 3.8& 4.5 &3.6 & 3.7 & 4.1 \\
        OR utilization \textsuperscript{1} (\%) & 85.5 & 85.2& 81.9&85.7 & 85.3 & 86.9\\
        \# Surgeries & 53 & 55& 51 & 54 & 53 & 55 \\
        Optimality gap (\%) &  0.71 & 0.95& 22.01 & 0.34 & 0.62 & 15.04 \\
        \bottomrule 
    \end{tabular} 
\end{subtable}

\vspace{1em} 

\begin{subtable}{\textwidth}
    \centering
    \small
    \begin{threeparttable}
    \caption{Performance ENT dataset 2} \label{tab:results_models_ENT2}
\begin{tabular}{p{3.5cm} p{1.3cm}p{1.3cm} p{1.3cm} p{1.3cm}p{1.3cm}p{1.3cm}}
        \toprule
        &  \multicolumn{3}{c}{\textbf{5 minutes}} &   \multicolumn{3}{c}{\textbf{1 hour}} \\
        \cmidrule(lr){2-4} \cmidrule(lr){5-7} 
        \textbf{KPI}  & FNN & PLF & SBM & FNN & PLF & SBM  \\ 
        \midrule 
        Objective value & 4,219 & 4,277& 4,211 &4,237 & 4,277 & 4,300\\
        \hangindent=1em \hangafter=0 Total duration: $\sum_{s \in \mathcal{S}}\bar{w}_s$ &  4,215 & 4,273& 4,207 &4,232 & 4,273 & 4,296\\
        \hspace{1em}Priority: $\sum_{s \in \mathcal{S}} \frac{1}{q_s+1}$ &  4.1 & 4.3& 4.3 &4.3 & 4.3 & 4.2 \\
        OR utilization \textsuperscript{1} (\%) & 79.4 & 80.5& 79.3 &79.7 & 80.5 & 80.9\\
        \# Surgeries &  37 & 41& 39 &40 & 41 & 39 \\
        Optimality gap (\%) &  2.09 & 1.22& 22.29 &1.09 & 1.21 & 17.38 \\
        \bottomrule 
    \end{tabular} 
    \begin{tablenotes}
\footnotesize
\item[] FNN: feedforward neural network; PLF: piecewise linear function; SBM: scenario-based modelling
\item[1] Based on the total mean scheduled duration and the total maximum OR capacity.
\end{tablenotes}
\end{threeparttable}
\end{subtable}
\end{table}

\subsection{Performance solution methods} \label{sec:comp_sol_methods}
\noindent
For each approach, we solve the models on the input datasets under time limits of both 5 minutes and 1 hour. Each run produces a schedule corresponding to the given dataset–approach combination. We evaluate these schedules in terms of their objective function value, split into total duration and priority, along with OR utilization, the number of scheduled surgeries, and optimality gap. 
Table \ref{tab:results_models} shows these metrics for all datasets.
The FNN and piecewise linear function approach show very similar results across all datasets and both time limits, with the FNN model yielding slightly higher objective function values in 6 of the 8 cases. 
When we look at scenario-based modelling, it yields a lower objective value for each dataset compared to the other approaches, except for both ENT datasets, under a time limit of 1 hour.
When zooming in on the total scheduled mean surgery duration of both cardiology datasets, scenario-based modelling allocates roughly 7 hours less for dataset 1 and 5 hours less for dataset 2 compared to both other approaches. 
This translates to almost a full dedicated OR less in capacity for dataset 1.
When we focus on the number of scheduled surgeries, we see that for the ENT datasets, each approach schedules roughly the same number of surgeries. For the cardiology datasets, however, we see that the scenario-based modelling approach schedules approximately 15-20 surgeries more than the other two approaches, while also showing a higher priority.
This indicates that incorporating uncertainty using different approaches leads to structurally different scheduling decisions.
Furthermore, we see that the optimality gap for the scenario-based modelling approach is higher than for the other approaches across all input datasets. 
This is due to the use of parameter $M_{od}^l$ in constraints \eqref{eq:constr_bigM}, which prevents the upper bound from decreasing, limiting the solver’s ability to tighten the gap. This leaves the solution with a weaker upper bound. The other two approaches are not restricted to such constraints, allowing their upper bounds to decrease and resulting in smaller optimality gaps.

In addition, Table \ref{tab:prob_overtime_constraints} reports the average overtime probabilities corresponding to the left-hand side of the overtime constraints \eqref{eq:fnn_overtime}, \eqref{eq:plf_overtime} and \eqref{eq:sbm_overtime} for each dataset–approach combination. The results indicate that the FNN approach yields overtime probabilities closest to the predefined threshold of $\alpha = 0.15$, followed by the piecewise linear function approach. In contrast, the scenario-based modelling approach consistently produces average overtime probabilities below 0.14 across all datasets and time limits. 
These findings suggest that, under the imposed time limits, the FNN approach utilizes available capacity most effectively while adhering to the approach’s overtime constraints.

\subsection{Simulation study}\label{sec:sim_study}
\noindent
Next, we conduct simulations to determine how often each generated schedule results in overtime based on realised durations. 
For these simulations, we use the general dataset from Section \ref{sec:generalprocess}.
For each generated schedule, we randomly select realised surgery durations from the general dataset corresponding to surgery type $t \in \mathcal{T}$ for each surgery $s \in \mathcal{S}_t$ on that schedule. 
The realised surgery duration $k_{od}$ for each operating room $o \in \mathcal{O}$ on day $d \in \mathcal{D}$ is then calculated as the sum of the selected realised durations. If the total realised surgery duration $k_{od}$ exceeds the maximum capacity $C_{od}$ of OR $o \in \mathcal{O}$ on day $d \in \mathcal{D}$, the generated schedule for that specific OR and day is recorded as an overtime occurrence.
We repeat this simulation 10,000 times for each generated schedule. The average overtime probabilities are computed and presented in Table \ref{tab:simulation}.
As visible, three of the generated schedules of scenario-based modelling result in an average overtime probability higher than 0.15, even though this approach has the lowest average overtime probabilities based on the approach's overtime constraints.
For the FNN approach, one of the generated schedules leads to an average overtime probability higher than 0.15. 
In addition, for the FNN and piecewise linear approach, we see that the average overtime probability imposed by the overtime constraints is higher than that of the simulation, except for the ENT 2 schedules. 
This indicates that, in general, these generated schedules are more conservative than implied by their overtime constraints.

Figure \ref{fig:simulation_nr_overtime_ORs} visualizes the number of ORs for each generated schedule that exhibit an overtime probability higher than 0.15 based on the simulation, denoted by excessive overtime. 
In general, all methods result in fewer operating rooms running in excessive overtime for the cardiology datasets than for the ENT datasets. 
This difference is likely explained by the lower mean standard deviation of the surgeries observed in the cardiology data relative to the ENT data.
Among the approaches, the piecewise linear function approach yields the smallest number of excessive overtime ORs, whereas the scenario-based modelling approach produces the highest number of ORs running in excessive overtime with even 80\% of ORs running in excessive overtime for ENT dataset 1 under a time limit of one hour.
This contrast further highlights the comparatively weaker overtime performance of the scenario-based modelling approach relative to the other approaches.

To evaluate the accuracy of each approach relative to real-life behaviour, we also compare the overtime probability derived from the overtime constraints with those obtained from the simulation on an OR level. This means that for each OR on each day, we compute the difference between the overtime probability based on the constraints and based on the simulation. Figure \ref{fig:diff_overtime_prob} shows the differences for each dataset-approach combination for both the schedules generated under a time limit of 5 minutes and 1 hour. A positive difference indicates that the overtime constraints show a higher probability than the simulation, while a negative difference indicates a higher simulated probability.

The first observation is that for most ORs, the overtime probability imposed by the constraints is higher than of the simulation.
This aligns with the average overtime probabilities shown in Table \ref{tab:probability_overtime_tabs}.
Next, we see that the difference for scenario-based modelling is larger than for the FNN and piecewise linear function approach. Those have quite similar differences in overtime probability across all datasets.
This shows that schedules generated under scenario-based modelling are less reliable than those produced by the FNN and piecewise-linear function approaches.
When we focus on the median for both cardiology and ENT datasets, we see comparable values for the FNN and piecewise linear function approach.
This indicates that despite OR-specific differences in overtime probability, the median of the difference is similar across the cardiology and ENT datasets for both approaches.
However, the median value under the scenario-based modelling approach deviates more from the other approaches, especially for both ENT datasets, which again shows that these schedules are less reliable.
Additionally, in the ENT 2 dataset, all generated schedules include one OR for which the simulated overtime probability is approximately 0.1 higher than the corresponding overtime-constraint probability. A closer inspection of the results shows that this arises from the same surgery being assigned to the operating room open from 08:00 to 20:00. This procedure has a mean duration of $\mu^{\text{(N)}} = 360$ and a standard deviation of $\sigma^{\text{(N)}} = 240.7$, which is the largest standard deviation observed across all datasets. 
This also explains the higher average overtime probabilities shown in Table \ref{tab:simulation}.
These findings suggest that all approaches encounter greater difficulty in accounting for surgeries with high duration variability.

\begin{table}[t]
\caption{Percentage overtime occurrences based on overtime constraints and simulation} 
\label{tab:probability_overtime_tabs}
\begin{subtable}{\textwidth}
    \centering
    \small
    \caption{Average overtime probability (\%) based on overtime constraints approaches} \label{tab:prob_overtime_constraints}
\begin{tabular}{l p{1.6cm} p{1.6cm} p{1.6cm} p{1.6cm} p{1.6cm} p{1.6cm}}
        \toprule
        &  \multicolumn{2}{c}{\textbf{FNN}} &   \multicolumn{2}{c}{\textbf{Piecewise linear function}} & \multicolumn{2}{c}{\textbf{Scenario-based modelling}} \\
        \cmidrule(lr){2-3} \cmidrule(lr){4-5} \cmidrule(lr){6-7}
        \textbf{Dataset}  & 5 min & 1 hour & 5 min & 1 hour & 5 min & 1 hour \\ 
        \midrule 
        Cardiology 1 &  14.9 & 14.9 & 14.5 & 14.6& 12.0 & 12.2 \\
        Cardiology 2 &  14.7 & 14.9 & 12.7 & 12.7& 11.2 & 11.3 \\
         ENT 1 &   14.9 & 14.9 & 14.2 & 14.3& 13.2 & 13.9 \\
         ENT 2 &     14.8 & 14.9 & 14.4 & 14.4& 13.6 & 13.7\\
        \bottomrule 
    \end{tabular} 
\end{subtable}
\vspace{0em} 

\begin{subtable}{\textwidth}
    \centering
    \small
    \caption{Average percentage (\%) of overtime occurrence based on simulation} \label{tab:simulation}
\begin{tabular}{l p{1.6cm} p{1.6cm} p{1.6cm} p{1.6cm} p{1.6cm} p{1.6cm}}
        \toprule
        &  \multicolumn{2}{c}{\textbf{FNN}} &   \multicolumn{2}{c}{\textbf{Piecewise linear function}} & \multicolumn{2}{c}{\textbf{Scenario-based modelling}} \\
        \cmidrule(lr){2-3} \cmidrule(lr){4-5} \cmidrule(lr){6-7}
        \textbf{Dataset}  & 5 min & 1 hour & 5 min & 1 hour & 5 min & 1 hour \\ 
        \midrule 
        Cardiology 1 &  12.6 & 13.1 & 12.7 & 12.3 & 10.1 & 10.7 \\
        Cardiology 2 &   12.5 & 12.8 & 10.9 &10.9 & 9.7 & 10.5 \\
         ENT 1 &  14.2 & 14.2 & 13.2 & 13.7 & 14.9 & \textcolor{red}{18.3} \\
         ENT 2 &   14.8 & \textcolor{red}{15.2} & 14.7 & 14.7 & \textcolor{red}{15.4} & \textcolor{red}{15.8}\\
        \bottomrule 
    \end{tabular} 
\end{subtable}
\end{table}

\begin{figure}[t]
\centering
\begin{tikzpicture}
\begin{axis}[
    ybar,
    ymin=0,
    width=16.9cm,
    height=6.5cm,
    bar width=12pt,
    enlarge x limits=0.09,
    symbolic x coords={
        A5, spacerA, A1, spacer1a, spacer1b, 
        B5, spacerB, B1, spacer2a, spacer2b, 
        C5,spacerC, C1, spacer3a, spacer3b, 
        D5, spacerD ,D1
    },
    xtick={A5,A1,B5,B1,C5,C1,D5,D1},
    xticklabels={
        \scriptsize 5 min, \scriptsize 1 hour,
        \scriptsize 5 min, \scriptsize 1 hour,
        \scriptsize 5 min, \scriptsize 1 hour,
        \scriptsize 5 min, \scriptsize 1 hour
    },
    ylabel={\small Nr ORs excessive overtime},
    xtick pos=bottom,
    ytick pos=left,
    legend style={
        font=\scriptsize,
        legend cell align=left,
        at={(0.99,0.98)},
        anchor=north east,
        draw=gray!70
    },
    extra x ticks={A5,B5,C5,D5},
    extra x tick labels={\footnotesize Cardiology 1,\footnotesize Cardiology 2,\footnotesize ENT 1,\footnotesize ENT 2},
    extra x tick style={
        tick label style={yshift=-12pt, xshift = 22pt, anchor=north},
        grid=none
    }
]

\addplot[
    fill=FNNcolor,
    bar shift=-12pt,
    nodes near coords,
    every node near coord/.style={font=\scriptsize, anchor=south, xshift=-12pt},
    point meta = explicit symbolic 
] coordinates {
    (A5,0)[] (A1,1)[1] (B5,1)[1] (B1,0)[] (C5,2)[2] (C1,2)[2] (D5,2)[2] (D1,3)[3]
};

\addplot[
    fill=PLFcolor,
    bar shift=0pt,
    nodes near coords,
    every node near coord/.style={font=\scriptsize, anchor=south},
    point meta = explicit symbolic 
] coordinates {
    (A5,1)[1] (A1,2)[2] (B5,0)[] (B1,0)[] (C5,1)[1] (C1,2)[2] (D5,1)[1] (D1,1)[1]
};

\addplot[
    fill=SBMcolor,
    bar shift=12pt,
    nodes near coords,
    every node near coord/.style={font=\scriptsize, anchor=south, xshift=12pt},
    point meta = explicit symbolic 
] coordinates {
    (A5,1)[1] (A1,1)[1] (B5,1)[1] (B1,4)[4] (C5,5)[5] (C1,8)[8] (D5,5)[5] (D1,4)[4]
};

\draw[dashed, gray, line width=0.3pt] (rel axis cs:0.25,0) -- (rel axis cs:0.25,1); 
\draw[dashed, gray, line width=0.3pt] (rel axis cs:0.5,0) -- (rel axis cs:0.5,1); 
\draw[dashed, gray, line width=0.3pt] (rel axis cs:0.75,0) -- (rel axis cs:0.75,1); 

\legend{
    FNN,
    Piecewise linear function,
    Scenario-based modelling
}

\end{axis}
\end{tikzpicture}
\caption{The number of ORs per dataset-approach combination that resulted in an overtime probability larger than 0.15 during the simulation. The cardiology setting has a total of 23 ORs, while the ENT setting contains 10 ORs.}
\label{fig:simulation_nr_overtime_ORs}
\end{figure}
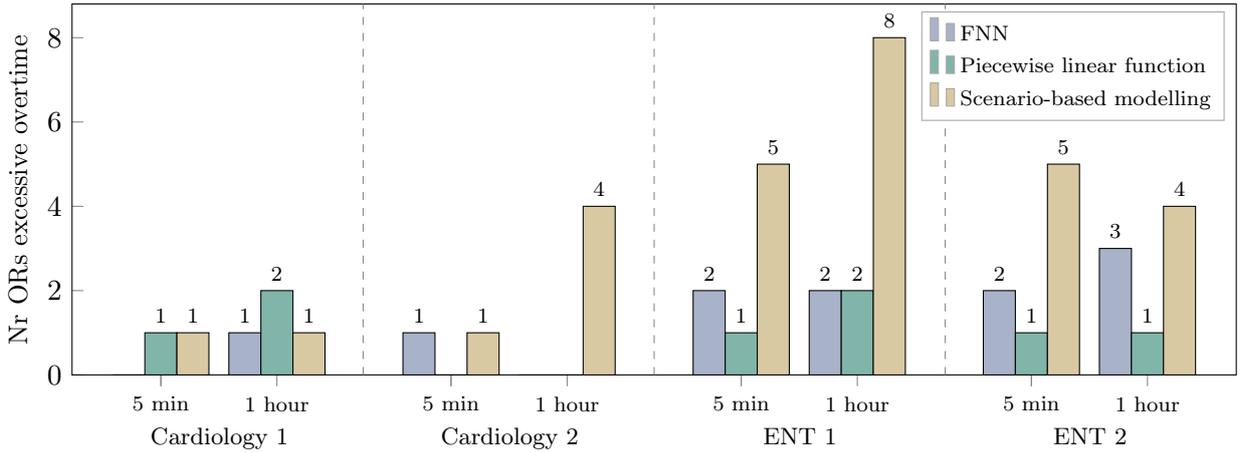

\begin{figure}[t]
\centering
\input{Figures/tikz_diff_prob_model_sim}
\caption{Difference of overtime probability overtime constraints - simulation on OR level. Cardiology settings contain 23 ORs, ENT settings contain 10 ORs.}
    \label{fig:diff_overtime_prob}
\end{figure}

\subsection{Feasibility generated schedules}\label{sec:feas_check}
\noindent
To determine the flexibility of the approaches, we evaluate whether the generated schedules under a certain approach are (partly) feasible solutions for both other approaches. 
We analyse the schedules that are generated under a time limit of 1 hour, since these yield the highest objective function values for all dataset-approach combinations.
Table \ref{tab:gen_schedules_feasibility_1hr} displays the feasibility of the solutions under the other approaches. The first column shows the approach with which the schedule is generated, while the other columns indicate the number of feasible ORs of the generated schedule under the other approaches.

\begin{table}[t]
\centering
\caption{Number of feasible OR schedules generated with time limit 1 hour for all approaches.}
\label{tab:gen_schedules_feasibility_1hr}
\begin{minipage}{0.45\textwidth}
\centering
\subcaption{Cardiology OR schedules (total ORs: 23)}

\begin{threeparttable}
\begin{tabular}{p{2.5cm} p{1.3cm}p{1.3cm} p{1.3cm} }
        \toprule

        Cardiology 1  & \textbf{FNN} & \textbf{PLF} & \textbf{SBM}   \\
        \midrule
        \hangindent=1em \hangafter=0 
        FNN & - & 11 & 9\\
        \hspace{1em}PLF & 23 & - & 8\\
        \hspace{1em}SBM & 12 & 12 & -\\
 
        \midrule
        Cardiology 2  &  \textbf{FNN} & \textbf{PLF} & \textbf{SBM}     \\
        \midrule
        \hangindent=1em \hangafter=0 
        FNN & - & 9 & 13\\
        \hspace{1em}PLF & 23 & - & 14\\
        \hspace{1em}SBM & 12 & 11 & -\\
    
        \midrule 
        \textbf{Total} & 70 & 43 & 44\\
        \bottomrule
    \end{tabular} 
\end{threeparttable}
\end{minipage}
\hfill
\begin{minipage}{0.45\textwidth}
\centering
\subcaption{ENT OR schedules (total ORs: 10)}
\begin{threeparttable}
\begin{tabular}{p{2.1cm} p{1.3cm}p{1.3cm} p{1.3cm} }
        \toprule
        ENT 1  &  \textbf{FNN} & \textbf{PLF} & \textbf{SBM}   \\
        \midrule
        \hangindent=1em \hangafter=0 
        FNN & - & 4 & 3\\
        \hspace{1em}PLF & 9 & - & 3\\
        \hspace{1em}SBM & 0 & 0 & -\\

        \midrule
        ENT 2  & \textbf{FNN} & \textbf{PLF} & \textbf{SBM}   \\
        \midrule
        \hangindent=1em \hangafter=0 
        FNN & - & 8 & 6\\
        \hspace{1em}PLF & 9 & - & 5\\
        \hspace{1em}SBM & 2 & 3 & -\\
    
        \midrule 
        \textbf{Total} & 20 & 15 & 17\\
        \bottomrule
    \end{tabular} 
\end{threeparttable}
\end{minipage}
\end{table}

For most generated schedules, at least one OR is classified as infeasible by the other approaches, with the exception of the two cardiology schedules produced using the piecewise linear function. These schedules are feasible solutions under the FNN approach. 
Furthermore, we see that none of the OR schedules generated by the scenario-based modelling approach for ENT 1 are feasible under the other two approaches. This is consistent with the simulation study, in which 8 of the 10 ORs resulted in excessive overtime, representing the dataset-approach combination with the most rejected ORs.
Lastly, we observe that overall, the FNN approach accepts the most OR schedules generated by the other approaches. 
Among the 90 OR schedules that were accepted by the FNN, 86 (96\%) were likewise accepted by the simulation. In comparison, of the 58 OR schedules accepted using the piecewise linear function approach, 55 (95\%) were also accepted by the simulation, while of the 61 OR schedules accepted through scenario-based modelling, 56 (92\%) were accepted based on the simulation.
This indicates that the FNN approach accepts a wider variety of solutions that are also deemed valid by (one of) the other two approaches as well as the simulation.

To illustrate the extent to which the other approaches accept or reject an OR schedule, Figure \ref{fig:feasibility_approaches} presents the overtime probabilities for each generated OR schedule under the other approaches.
In general, we see that the OR schedules generated by the scenario-based model lead to the highest overtime probabilities under the other approaches, visualized in Figure \ref{fig:feasibility_SBM}. 
The scenario-based modelling approach also leads to a wider range of overtime probabilities for the schedules generated under the other approaches, indicating that this approach exhibits higher variability in both generating and accepting schedules.
Furthermore, the OR schedules generated by the FNN and piecewise linear function approaches that are rejected by either method, only exhibit overtime probabilities slightly above 0.15, where the OR schedules generated by the piecewise linear function method show the smallest overtime probabilities under the FNN approach.
This suggests that the schedules generated by the piecewise linear function method are closer to the FNN feasibility threshold than the reverse case, which aligns with the higher number of OR schedules accepted by the FNN approach.

Lastly, a noteworthy point is that even though scenario-based modelling assigns approximately 7 hours and 5 hours less surgery time for the cardiology datasets, only half of the ORs in the generated schedules are feasible for both other approaches. 
This once again shows that having fewer scheduled surgery hours does not guarantee a more reliable schedule.

\begin{figure}[!t]
    \centering
    \input{Figures/tikz_feasibility_approaches}
    \caption{Overtime probability for each OR and day of schedules generated under 1 hour for other approaches. Horizontal line indicates $\alpha = 0.15$.}
    \label{fig:feasibility_approaches}
\end{figure}

\subsection{Comparing solution methods with real-life overtime cases}\label{sec:real-life-cases}
\noindent
To determine if our approaches accommodate certain overtime cases with an overtime probability smaller than 0.15, we analyse cases in which the hospital’s initial schedule on OR level resulted in overtime.
From the general dataset introduced in Section \ref{sec:generalprocess}, we identify all realised overtime cases for which more than one surgery was scheduled in the same operating room. 
An overtime case is defined as a case where the realised total surgery duration exceeds the OR capacity on a given day.
This leads to 160 overtime cases for cardiology and 49 overtime cases for ENT. The maximum OR capacity is 510 minutes for all these cases.

We divide the overtime cases into two groups based on a simulation: cases with an average simulated overtime probability of at most 0.15 and cases with an average probability exceeding 0.15. The simulation procedure follows the approach described in Section \ref{sec:comp_sol_methods}, using a total of 10,000 simulations for each overtime case. 
Among the 160 cardiology cases, 37 have an average overtime probability of at most 0.15, whereas this holds for 9 of the 49 ENT cases.

To assess whether the approaches accept cases with a simulated average overtime probability of at most 0.15 while rejecting the remaining cases, we evaluate each overtime case individually.
For all runs, consider one available OR, a one-day planning horizon, and a maximum capacity $C_{od}$ of 510 minutes.
Furthermore, surgery set $\mathcal{S}$ consists of all surgeries that were scheduled in a specific overtime case. For each of these surgeries $s \in \mathcal{S}$, we set release date $r_s$ equal to 0, due date $q_s$ equal to 1 and priority $p_s$ equal to 0.
We apply the objective function presented in equation \eqref{eq:generalobjfunc}. 
If all surgeries in set $\mathcal{S}$ are scheduled according to the capacity constraints, the approach allows the overtime case under overtime probability $\alpha = 0.15$.
If only a subset of surgeries is scheduled, the approach does not allow the specific overtime case.

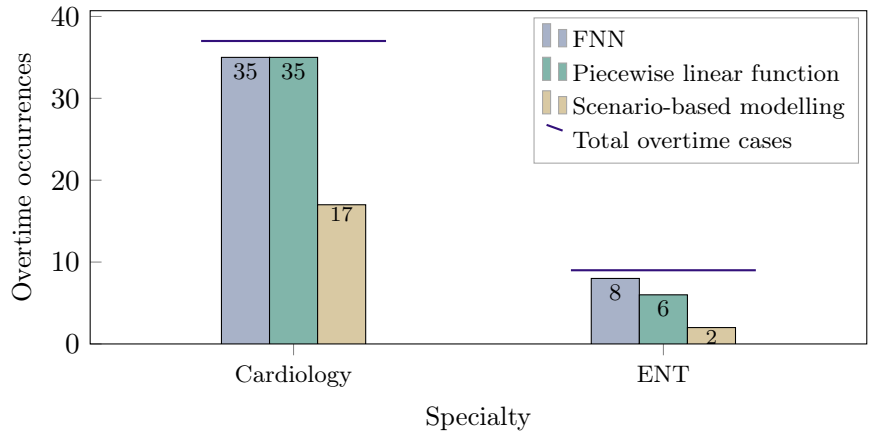
\begin{figure}[t]
\centering
\begin{tikzpicture}
\begin{axis}[
    ybar,
    ymin=0,
    width=11.8cm,
    height=6cm,
    bar width=18pt,
    enlarge x limits=0.2,
    xtick={0,1},
    xticklabels={\footnotesize Cardiology, \footnotesize ENT},
    legend style={
        font=\footnotesize,
        legend cell align=left,
        at={(0.99,0.98)},
        anchor=north east,
        draw=gray!70
    },
    ylabel={\small Overtime occurrences},
    xlabel={\small Specialty},
    xtick pos=bottom,
    ytick pos=left
]

\addplot[
    fill=FNNcolor,
    bar shift=-18pt,
    nodes near coords,
    every node near coord/.style={
        font=\footnotesize,
        anchor=south,
        xshift=-18pt,   
        yshift = -12pt
    }
] coordinates {(0,35) (1,8)};

\addplot[
    fill=PLFcolor,
    bar shift=0pt,
    nodes near coords,
    every node near coord/.style={
        font=\footnotesize,
        anchor=south,
        yshift = -12pt
    }
] coordinates {(0,35) (1,6)};

\addplot[
    fill=SBMcolor,
    bar shift=18pt,
    nodes near coords,
    every node near coord/.style={
        font=\scriptsize,
        anchor=south,
        xshift=18pt,   
        yshift = -10pt
    }
] coordinates {(0,17) (1,2)};

\addplot[
    sharp plot,
    color=CAR1color,
    thick
] coordinates {
    (-0.25,37)
    (0.25,37)
};

\addplot[
    sharp plot,
    color=CAR1color,
    thick
] coordinates {
    (0.75,9)
    (1.25,9)
};

\legend{
    FNN,
    Piecewise linear function,
    Scenario-based modelling,
    Total overtime cases
}

\end{axis}
\end{tikzpicture}
\caption{The number of accepted overtime cases per approach versus the total number of overtime cases with a simulated overtime probability of at most 0.15.}
\label{fig:comp_overtime_cases_alph015}
\end{figure}

From the overtime cases for which the simulated average overtime probability is greater than 0.15, none were accepted by any of the approaches.
This indicates that all approaches avoid accepting these overtime cases with unacceptably high overtime risk.
When zooming in on the overtime cases with a simulated average overtime probability of at most 0.15, we see that each approach performs differently.
Figure \ref{fig:comp_overtime_cases_alph015} shows the number of overtime cases that are allowed by an approach, where the horizontal line indicates the total number of these cases. The scenario-based modelling approach allows the fewest cases, accepting a total of 19 overtime cases. In comparison, the piecewise linear function approach allows 41 cases, while the FNN approach allows 43 cases. 

These results demonstrate that the scenario-based modelling approach is highly conservative for these cases, as evidenced by the substantial gap between the number of cases with an average overtime probability of at most 0.15 and the number of cases it allows. The FNN and piecewise linear function approaches yield similar results and are closer to this benchmark, with the FNN approach performing slightly better.
Thus, the FNN method allows the most cases that are also deemed acceptable by the simulation.

\section{Discussion and conclusion}\label{Sec:discussion}
\noindent
In this paper, we analyse and compare three different approaches to include uncertainty of surgery durations within MILP models.
All methods provide a framework to incorporate variability in surgery durations, but differ in their performance under uncertainty and level of conservatism, which affects both overtime probabilities and OR utilization. 

We first focus on the FNN approach, which uses the Fenton-Wilkinson method to approximate the sum of lognormal distributions. This method performs better for datasets with low variability, as the approximation becomes less accurate when the summed variables present high variances. 
This may explain why the number of ORs that experience a simulated overtime probability larger than 0.15 is higher for the ENT datasets than for the cardiology datasets, given that the surgery durations in the ENT datasets show greater variability.
However, when we focus on the number of surgery types that better fit a lognormal distribution versus a normal distribution, the ENT datasets show a better fit. This could explain why the difference in overtime probability between constraints and simulation is smaller for the ENT datasets compared to the cardiology datasets.

Examining the characteristics of the lognormal distribution further, we observe that when the lognormal standard deviation parameter decreases, the distribution becomes more symmetric and increasingly concentrated around its mean. As skewness decreases, the lognormal distribution more closely resembles a normal distribution.
Translating this insight to the piecewise linear function and FNN approach, we would expect the piecewise linear method to yield simulation results that are closer to those of the FNN method for the cardiology datasets than for the ENT datasets.
This pattern is indeed observed in the cardiology datasets, where the median difference between constrained and simulated overtime probabilities for the piecewise linear approach is closer to zero compared to the ENT datasets. However, for certain individual ORs, the piecewise linear approach still exhibits larger deviations. When we focus on the ENT datasets, both approaches show small differences.

When we compare the FNN approach with the piecewise linear function model on objective function value, we see that in 6 of the 8 generated schedules, the objective values obtained by the FNN approach are slightly higher. 
The FNN and piecewise linear function approach have similar total number of overtime ORs based on simulations, where the FNN outperforms in the cardiology setting and the piecewise linear function approach show better results in the ENT setting.
In contrast, the scenario-based modelling approach is generally outperformed by the other two approaches with respect to both objective function value and simulated overtime probabilities. Although scenario-based modelling achieves the lowest objective function value in the cardiology 2 dataset, it results in the largest number of ORs with a simulated overtime probability above 0.15. A similar pattern is observed in the ENT 1 dataset, where a substantial proportion of generated OR schedules exceed the 0.15 overtime threshold in simulation, despite the average overtime probability implied by the scenario-based model’s constraints being the lowest among the three approaches. 
So, in general, although scenario-based modelling results in lower objective function values, it also produces the highest number of overtime ORs in simulation. 
Even if the method were allowed a longer runtime to potentially improve the objective value, the approach already results in considerable OR overtime.
These findings indicate that the scenario-based modelling approach may not translate as effectively to simulated real-world performance and appears less reliable in controlling overtime risk.

Another important aspect concerns the search space and feasibility of the generated schedules.
The FNN approach accepts a wider variety of solutions proposed by the other approaches that are also deemed acceptable based on simulation. 
As a result, it is able to explore a wider range of feasible solutions to maximize utilization under the given constraints. In contrast, schedules produced by the scenario-based modelling approach are rejected most frequently by the other two approaches as well as by the simulation study. This suggests that scenario-based modelling tends to generate schedules that are less often considered feasible under the criteria applied by the other methods and the simulation study.

Under all approaches, real-life overtime cases for which the simulated overtime probability exceeds 0.15 are rejected. This indicates that, for these cases, the approaches adhere more strictly to the predefined probability threshold than the schedules observed in practice.
Furthermore, when we focus on the overtime cases that were allowed based on the simulated overtime probability, scenario-based modelling appears the most conservative, as it permits the fewest of these overtime cases.
The FNN approach accepts the most of these overtime cases with 2 more allowed cases than the piecewise linear function method.
This again shows that the FNN model allows a wider set of options that are also deemed feasible based on simulations.

Taken together, the results suggest that the FNN and piecewise linear function approaches produce more reliable and practically applicable schedules than the traditional scenario-based modelling method. Within the evaluated context, the FNN approach appears most suitable for incorporating uncertainty in surgery duration within MILP models. It achieves the highest objective function values, maintains competitive simulated overtime probabilities relative to the piecewise linear approach, and accepts the largest number of schedules proposed by alternative methods. Moreover, it yields average simulated overtime probabilities that are closest to the predefined threshold of 0.15 and accommodates the greatest number of real-life overtime cases that remain below this threshold in simulation.
A limitation is that the FNN approach, like the other two methods, still generates some OR schedules with simulated overtime probabilities exceeding 0.15. However, this limitation is not unique to FNN. 
In general, the average simulated overtime probability imposed by the FNN approach is closest to the predefined overtime probability for all datasets.
This shows that this approach has the smallest average deviation from this predefined target compared to the other approaches.
Overall, the FNN approach appears to offer the most balanced trade-off between utilization and overtime risk under the given modelling framework, while the piecewise linear function approach remains a viable and competitive alternative depending on dataset characteristics.

\section*{Declaration of interest}
\noindent
Declarations of interest: none


\section*{Declaration of generative AI and AI-assisted technologies}
\noindent
During the preparation of this work, the authors used ChatGPT in order to refine the wording and improve the clarity and phrasing of the manuscript. The tool was employed for the lay-out of figures and language editing purposes, including enhancing readability and ensuring consistency in expression. Furthermore, ChatGPT was used to debug parts of the code more easily. After using this tool, the authors reviewed and edited the content as needed and take full responsibility for the content of the published article.

\section*{Funding}
\noindent
This research is supported by Convergence Health and Technology Sustainable Health Program SMART OR2030 [grant number 2023014].

\printbibliography
\end{document}